\documentclass[11pt]{article}
\usepackage{amsmath}
\usepackage{amssymb}
\usepackage{amsbsy}
\usepackage{amsthm}
\usepackage{epsfig}
\usepackage{wrapfig}
\usepackage{eepic}
\usepackage{color}
\usepackage{graphicx}
\usepackage{amsfonts}%

\newtheorem{Theorem}{Theorem}[section]

\newtheorem{remark}{Remark}[section]
\newcommand{\bu}{{\bf u}}

\newcommand{\be}{{\bf e}}
\newcommand{\bD}{{\bf D}}

\newcommand{\bv}{{\bf v}}
\newcommand{\bw}{{\bf w}}
\newcommand{\bx}{{\bf x}}
\newcommand{\bn}{{\bf n}}

\newcommand{\bH}{{\bf H}}
\newcommand{\bX}{{\bf X}}
\newcommand{\bL}{{\bf L}}
\newcommand{\bV}{{\bf V}}

\newcommand{\bpsi}{{\boldsymbol \psi}}
\newcommand{\bfeta}{{\boldsymbol \eta}}

\newcommand{\bphi}{{\boldsymbol \phi}}

\def\div{\operatorname{div}}

\newcommand{\blf} {\mathbf{f}}

\newcommand{\MO}[1]{{\color{black}#1}}

\begin{document}

\title{Longer time accuracy for incompressible Navier-Stokes simulations with the EMAC formulation}

\author{
Maxim A. Olshanskii
\footnote{Department of Mathematics, University of Houston, Houston TX 77004 (molshan@math.uh.edu)}
\and
Leo G. Rebholz\footnote{Department of Mathematical Sciences, Clemson University, Clemson, SC 29634 (rebholz@clemson.edu)}
}
\date{}

\maketitle

\begin{abstract}
In this paper, we consider the recently introduced  EMAC formulation for the incompressible Navier-Stokes (NS) equations, which is the only known NS formulation that conserves
energy, momentum and angular momentum when the divergence constraint is only weakly enforced.  Since its introduction, the EMAC formulation  has been successfully used for a wide variety of fluid dynamics problems.  We prove that discretizations using the EMAC formulation are potentially better than those built on the commonly used skew-symmetric formulation, by deriving a better longer time error estimate for EMAC: while the classical results for schemes using the skew-symmetric formulation have Gronwall constants dependent on $\exp(C\cdot Re\cdot T)$ with $Re$ the Reynolds number, it turns out that the EMAC error estimate is free from this explicit exponential dependence on the Reynolds number.  Additionally, it is demonstrated how EMAC admits  smaller lower bounds on  its velocity error, since {incorrect treatment of linear momentum, angular momentum and energy induces} lower bounds for $L^2$ velocity error, and EMAC treats these quantities more accurately.  Results of numerical tests for channel flow past a cylinder and 2D Kelvin-Helmholtz instability are also given, both of which show that the advantages of EMAC over  the  skew-symmetric formulation increase as the Reynolds number gets larger and for longer simulation times.
\end{abstract}

\section{Introduction}

We consider herein numerical schemes for the incompressible Navier-Stokes equations (NSE), which are given by
\begin{eqnarray}
\bu_t + (\bu\cdot\nabla) \bu + \nabla p - \nu\Delta \bu & = & \blf, \label{nse1} \\
\div\bu & = & 0, \label{nse2}\\
\bu(0) &= &\bu_0,
\end{eqnarray}
in a domain $\Omega\subset \mathbb{R}^d$, $d$=2 or 3, with polyhedral and Lipschitz boundary, $\bu$ and $p$ representing the unknown velocity and pressure, $\blf$ an external force, $\bu_0$ the initial velocity, and $\nu$ the kinematic viscosity which is inversely proportional to the Reynolds number $Re$.  Appropriate boundary conditions are required to close the system, and for simplicity we will consider the case
of homogeneous Dirichlet boundary conditions, $\bu |_{\partial\Omega}={\bf 0}$.

In the recent work \cite{CHOR17}, the authors showed that due to the divergence constraint, the NSE nonlinearity could be equivalently be written as
\[
\bu\cdot\nabla \bu + \nabla p = 2\bD(\bu)\bu + (\div \bu)\bu + \nabla P,
\]
with $P=p-\frac12 |\bu|^2$ and $\bD$ denoting the rate of deformation tensor.  Reformulating in this way was named
in \cite{CHOR17} to be the {\it energy, momentum and angular momentum conserving} (EMAC) formulation of the NSE, since when discretized with a Galerkin method that only weakly enforces the divergence constraint, the EMAC formulation still produces a scheme that conserves each of energy, momentum, and angular-momentum, as well as properly defined 2D enstrophy, helicity, and total vorticity.  This is in contrast to the well-known convective, conservative, rotational, and skew-symmetric formulations, which are each shown in \cite{CHOR17} to not conserve at least one of energy, momentum or angular momentum.

The EMAC formulation, and its related numerical schemes, is part of a long line of research extending back at least to Arakawa that has the theme ``incorporating more accurate physics into discretizations leads to more stable and accurate numerical solutions, especially over long time intervals.''  There are typically many ways to discretize a particular PDE, but choosing (or developing) a method that more accurately reproduces important physical balances or conservation laws will often lead to better solutions.  Arakawa recognized this when he designed an energy and enstrophy conserving scheme for the 2D Navier-Stokes equations in \cite{A66}, as did Fix for ocean circulation models in \cite{F75}, Arakawa and Lamb for the shallow water equations \cite{AL81}, and many others {for various evolutionary systems from physics}, e.g. \cite{ST89,AM03,LW04,R07,evans2013isogeometric,SCN15,PG16}.  It is important to note that if pointwise divergence-free approximations are used, such as those recently developed in \cite{arnold:qin:scott:vogelius:2D,Z05,evans2013isogeometric_s,GN14b,GN14}, then the numerical velocity found with the EMAC formulation is the same vector field as recovered from more traditional convective and skew-symmetric formulations, and all of these conservation properties will hold for those formulations as well.  However, the development of strongly divergence-free methods is still quite new, often requires non-standard meshing and elements, and is not yet included into most major software packages.

Since its original development in 2017 in \cite{CHOR17}, the EMAC formulation has gained considerable attention by the CFD community.  It has been used for a wide variety of problems, including vortex-induced vibration \cite{PCLRH18}, turbulent flow simulation \cite{LHOCR19}, cardiovascular simulations and hemodynamics \cite{SP18,SP18b}, noise radiated by an open cavity \cite{MSLGD19}, and others \cite{OCAM19,LPH19}.  It has proven successful in these simulations, and a common theme reported for it has been that it exhibits low dissipation compared to other common schemes, which is likely due to EMAC's better adherence to physical conservation laws and balances.  Not surprisingly, less has been done from an analysis viewpoint, as only one paper has appeared in this direction; in \cite{CHOR19}, the authors analyzed conservation properties of various time stepping methods for EMAC. In particular,  no analysis for EMAC has been found which improves upon the well-known analysis of mixed finite elements for the incompressible NSE in skew-symmetric form.  The present paper addresses the challenge of providing such new analysis.

This paper extends the study of the EMAC formulation  both analytically and computationally.  Analytically, we show how the better convergence properties of EMAC unlock the potential for decreasing the approximation error of FE methods. In particular, we show that while the classical semidiscrete error bound for the skew-symmetric formulation has a Gronwall constant $\exp(C \cdot Re \cdot T)$ \cite{J16}, where $T$ is the {simulation} end time, the analogous EMAC scheme has a Gronwall constant $\exp(C \cdot T)$, i.e. with no explicit exponential dependence on $Re$ (and the rest of the terms in the error bound are similar).  We note that previously, such $\nu$-uniform error bounds were believed to be an exclusive property of {finite element methods that enforced the divergence constraint strongly through divergence-free elements \cite{SLLL18} or through stabilization/penalization of the
divergence error \cite{FGJN18}.}  Additionally, we show how the lack of momentum conservation in convective, skew-symmetric and rotational forms produce a lower bound on the error, which EMAC is free from.  Numerically, we first extend results from \cite{CHOR17} for flow past a cylinder to higher Reynolds number, and also consider Kelvin-Helmholtz instability simulations.  For both tests, we find that EMAC produces better results than the analogous skew-symmetric scheme when the flow is underresolved - which is the case of practical interest.

This paper is arranged as follows. In section 2, we provide notation and mathematical preliminaries for a smooth analysis to follow.  Section 3 presents new analytical results for EMAC, including proving an improved Gronwall constant that is not explicitly dependent on the Reynolds number and showing that a lack of conservation laws (i.e. for non-EMAC schemes) creates a lower bound on $L^2$ velocity error.  Section 4 presents results for two challenging numerical tests, both of which reveal advantages of EMAC for higher $Re$ problems.  Conclusions and future directions are discussed in section 5.

\section{Notation and preliminaries}

We consider a convex  polygon or polyhedral domain $\Omega \subset \mathbb{R}^d$, $d=2,3$.  We denote the $L^2(\Omega)$ inner product and norm on $\Omega$ by $(\cdot,\cdot)$ and $\| \cdot \|$, respectively.  The natural velocity and pressure spaces are
\begin{eqnarray*}
\bX =  \{ \bv \in H^1(\Omega)^d,\ \bv |_{\partial\Omega}= {\bf 0} \}, \qquad
Q =  \{ q \in L^2(\Omega),\ \int_{\Omega} q\ dx=0 \}.
\end{eqnarray*}
Let $\bV$ denote the divergence-free subspace of $\bX$,
$
\bV:= \{ \bv \in \bX\,:\,\nabla \cdot \bv=0, \textrm{a.e. in}~\Omega\}.
$

We will consider subspaces $\bX_h \subset \bX$, $Q_h \subset Q$ to be finite element (FE) velocity and pressure spaces corresponding to an admissible triangulation $\mathcal{T}_h$ of $\Omega$, where $h$ is the global mesh-size. For $\mathcal{T}_h$ we assume the minimal angle condition if $h$ varies.  We further assume that $\bX_h$ and $Q_h$ satisfy the inf-sup compatibility condition~\cite{GR86} and define the discretely divergence-free subspace of $\bX_h$ by
\[
\bV_h:= \{ \bv_h \in \bX_h,\ (\nabla \cdot \bv_h,q_h)=0\ \forall q_h\in Q_h \}.
\]
Most common FE discretizations of the NSE and related systems, e.g. using Taylor-Hood elements, only enforce the divergence constraint $\div \bu_h=0$ weakly and thus $\bV_h \not\subset \bV$.

We denote by $I_{St}^h$ the discrete Stokes projection operator, which is defined by: Given $\bw\in \bV$, find $I_h^{St}\bw \in \bV_h$ satisfying
\[
(\nabla I_h^{St}\bw,\nabla \bv_h) = (\nabla \bw,\nabla \bv_h) \quad \forall \bv_h \in \bV_h.
\]
The optimal order approximation properties of $I_h^{St}$ in $L^2$ and $H^1$ norms  readily follows from the inf-sup compatibility and interpolation properties of the finite element spaces~\cite{GR86}.
We further need stability of the Stokes projection in $W^{1,r}$ norms:
\begin{equation}
\| \nabla I_h^{St} \bw \|_{L^{r}} \le C \| \nabla \bw \|_{L^{r}},\quad r\in[2,\infty].  \label{stokesbound}
\end{equation}
This estimate is shown in Theorem 13 and Corollary 4 of \cite{girault2015max} under the above assumptions on  $\Omega$ and $\mathcal{T}_h$ if the pair $\bV_h,Q_h$ admits the construction of the Fortin projector with a local stability property. In turn, the existence of such a projector was demonstrated, cf. \cite{girault2015max}, for many popular inf-sup stable elements, including the Taylor--Hood element $P_{k}-P_{k-1}$ {for $k\ge d$}, {\color{black} where $P_k$ represents piecewise continuous polynomials of degree $k$ on each element},
and for the $k=2$, $d=3$ if the triangulation $\mathcal{T}_h$ has a certain macrostructure.

\subsection{Vector identities and trilinear forms}
For a sufficiently smooth velocity field $\bu$, the symmetric part of its gradient, $\bD(\bu) = \frac12 \left(\nabla \bu + (\nabla \bu)^T\right)$, defines the rate of deformation tensor. The EMAC formulation employs the following identity
\begin{equation}
(\bu\cdot \nabla) \bu = 2\bD(\bu)\bu - \frac12\nabla |\bu|^2, \label{vecid7}
\end{equation}
which splits the inertia term into the acceleration driven by $2\bD(\bu)$ and potential term further absorbed by redefined pressure.
Based on \eqref{vecid7} one defines the trilinear form for EMAC Galerkin formulation:
\[
c(\bu,\bv,\bw) =2(\bD(\bu)\bv,\bw)+(\mbox{div}(\bu)\bv,\bw),
\]
here and further in this section $\bu,\ \bv,\ \bw\in \bX$ (no divergence free condition is assumed for any of vector fields).
The divergence term is added in the definition of $c(\cdot,\cdot,\cdot)$ to ensure the cancellation property:
\begin{equation}\label{skew}
c(\bv,\bv,\bv)=0, 
\end{equation}
which leads energy balance without the div-free condition strongly enforced.
The popular convective form of the Galerkin method uses
\[
b(\bu,\bv,\bw) = (\bu\cdot\nabla \bv,\bw)
\]
for the nonlinear part of the equations, while skew-symmetric form anti-symmetrizes $b$,
\[
b^*(\bu,\bv,\bw)  = \frac12 ( b(\bu,\bv,\bw) - b(\bu,\bw,\bv)),
\]
to ensure the energy conservation through $b^*(\bu,\bv,\bv)=0$. The EMAC form of nonlinear terms, however, {also conserves} linear and angular momenta; see~\cite{CHOR17}.

We are now prepared to introduce finite element formulations.

\subsection{Semi-discrete FEM formulations}

We shall consider the following semi-discretization of the NSE, which uses the EMAC formulation of the nonlinear term:  Find $(\bu_h,P_h)\in (\bX_h,Q_h) \times (0,T]$ satisfying for all $(\bv_h,q_h)\in (\bX_h,Q_h)$,
\begin{align}
(({\bu_h})_t,\bv_h) + c(\bu_h,\bu_h,\bv_h) - (P_h,\nabla \cdot \bv_h) + \nu (\nabla \bu_h,\nabla \bv_h) & = (\blf,\bv_h), \label{emac1} \\
(\nabla \cdot \bu_h,q_h) &=0. \label{emac2}
\end{align}
We refer to \eqref{emac1}-\eqref{emac2} as the EMAC scheme.

In order to compare EMAC results, we shall also consider in our analysis the analogous  scheme with skew-symmetric form of nonlinear terms,  and will refer to it as SKEW.  The SKEW formulation is the same as \eqref{emac1}-\eqref{emac2}, except replacing  $c(\bu_h,\bu_h,\bv_h)$ with $ b^*(\bu_h,\bu_h,\bv_h)$, and renaming the pressure $p_h$ since it now represents the kinematic pressure:
Find $(\bu_h,p_h)\in (\bX_h,Q_h) \times (0,T]$ satisfying for all $(\bv_h,q_h)\in (\bX_h,Q_h)$,
\begin{align}
(({\bu_h})_t,\bv_h) + b^*(\bu_h,\bu_h,\bv_h) - (p_h,\nabla \cdot \bv_h) + \nu (\nabla \bu_h,\nabla \bv_h) & = (\blf,\bv_h), \label{skew1} \\
(\nabla \cdot \bu_h,q_h) &=0. \label{skew2}
\end{align}
The convective scheme (CONV) is precisely {\color{black} SKEW but without the skew-symmetry, i.e. \eqref{skew1}-\eqref{skew2} with
$b^*(\bu_h,\bu_h,\bv_h)$ replaced by $b(\bu_h,\bu_h,\bv_h)$. The rotational scheme (ROT) is created by replacing   $c(\bu_h,\bu_h,\bv_h)$ with $ ((\nabla \times \bu_h)\times \bu_h, \bv_h)$ in  \eqref{skew1}-\eqref{skew2} but here the pressure $P_h$ represents Bernoulli pressure, and the conservative scheme (CONS) is created by replacing   $b^*(\bu_h,\bu_h,\bv_h)$ with $b(\bu_h,\bu_h,\bv_h)+((\div \bu_h)\bu_h,\bv_h)$  in  \eqref{skew1}-\eqref{skew2}; see also Table~\ref{tab1}.}

Numerical tests from \cite{CHOR17,CHOR19} show that CONV and CONS, neither of which conserve energy (see section 2.3 below), can become unstable in problems with higher Reynolds numbers at longer time intervals, at least if no additional stabilization terms are added to the formulation.  Furthermore, ROT is well known to be less accurate in many cases compared to the other schemes \cite{LMNOR09} if common element choices are made.  For these reasons, this paper restricts analytical and numerical comparisons of EMAC only to SKEW.

\subsection{Energy balance and conservation of linear and angular  momentum by NSE}

Smooth solutions to NSE are well known  to deliver energy balance and conserve  linear momentum, and angular momentum, which we define by
\begin{eqnarray*}
\text{Kinetic energy}\quad && E:= \frac12 \int_\Omega |\bu|^2\mbox{d}\bx;\\
\text{Linear momentum}\quad&& M:=\int_\Omega \bu\, \mbox{d}\bx;\\
\text{Angular momentum}\quad&& M_\bx:=\int_\Omega \bu\times\bx\, \mbox{d}\bx.
\end{eqnarray*}
To see the balances, assume for simplicity that the solution $\bu,\ p$ have compact support in $\Omega$ (e.g. consider an isolated vortex), and test the NSE with $\bu$, {$\bpsi_i$ (the $i^{th}$ standard basis vector), and $\bphi_i =\bpsi_i \times \bx$} to obtain
\[
\frac{d}{dt}E + \nu \| \nabla \bu \|^2 =  (\blf,\bu), \quad
\frac{d}{dt} M_i =  \int_{\Omega} f_i \ d\bx, \quad
\frac{d}{dt} (M_\bx)_i = \int_{\Omega} (\blf\times \bx)_i\ d\bx,
\]
noting that each nonlinear and pressure term vanished, and using
\[
(\bu_t,\bphi_i)=\frac{d}{dt} \int_{\Omega} (\bu \times \bx)\cdot {\bpsi_i} \ d\bx = \frac{d}{dt} (M_\bx)_i.
\]
For a numerical scheme to have physical accuracy, its solutions should admit balances that match these balances as close as
possible.  The key point here is that the nonlinear term does not contribute to any of these balances.

It is shown in \cite{CHOR17} that EMAC conserves each of energy, momentum, and angular momentum when the divergence constraint is only weakly enforced, while none of the other schemes do (unless $\nabla \cdot \bu_h=0$ is strongly enforced, in which case all schemes are equivalent and all conserve each of these quantities).  Table \ref{tab1}, which is given in \cite{CHOR17}, summarizes the conservation properties of these schemes as well as their nonlinearity form and potential (pressure) term.

\begin{table}[h!]

\begin{tabular}{rl|c|c|c|c}
name& nonlinear form  & potential term                                  &  Energy & Mom. & Ang. Mom. \\ \hline
CONV: & $\bu \cdot\nabla \bu$ &$p$ (kinematic)             &&&\\
SKEW: & $\bu \cdot\nabla \bu +  \frac12 (\div \bu)\bu$&$p$ (kinematic) &\textbf{+}&&\\
ROT: & $(\nabla\times\bu)\times \bu$&$p+\frac12|\bu|^2$ (Bernoulli)       &\textbf{+}&&\\
CONS: & $\nabla\cdot (\bu{\color{black}\otimes}\bu)$&$p$ (kinematic)&&\textbf{+}&\textbf{+}\\
EMAC: &$2D(\bu)\bu + (\div \bu)\bu$&$p-\frac12|\bu|^2$ (no name)  &\textbf{+}&\textbf{+}&\textbf{+}
\end{tabular}
\caption{\label{tab1} Properties of the common NSE schemes}
\end{table}

\section{Error analysis}

We show in this section two improvements EMAC provides over SKEW in terms of error analysis. Putting our effort into a more general context, the paper addresses the challenge of providing mathematical evidence of the widely accepted concept that {\it a better adherence to conservation laws with discretization schemes leads to more accurate numerical solutions}.  Numerical tests using EMAC in, e.g., \cite{CHOR17,CHOR19,PCLRH18,LHOCR19,OCAM19,SP18,SP18b,LPH19,MSLGD19} already show it performs very well on a wide variety of application problems, and in many instances better than other formulations, as it is found to be less numerically dissipative; hence this section gives some analytical backing to those results.

\subsection{A lower bound for velocity error based on violating momentum/energy conservation}

Here we present a discussion of how a violation of conservation laws yields a lower bound on velocity error.  This is
a positive result for EMAC in the sense that EMAC is expected to have smaller momentum, angular momentum and energy errors than SKEW, ROT, CONS and CONV, since
the EMAC nonlinearity preserves all of these quantities while the other formulations violate the conservation of at least one of them.

For momentum conservation, the argument works as follows: let $\bu$ be the true solution, $\bu_h$ the discrete solution, and $\bpsi_i$ be the $i^{th}$ standard basis vector.  Then denoting
$e_i^{mom}$ to be the momentum error in the $i^{th}$ component,
\[
e^{mom}_i(t) = \bigg| \int_{\Omega} (\bu(t) - \bu_h(t))\cdot \bpsi_i\ d\bx \bigg| \le \int_{\Omega} \bigg|  (\bu(t) - \bu_h(t)) \bigg|  \ d\bx = \| \bu(t) - \bu_h(t) \|_{L^1}.
\]
Since $\|  \bu(t) - \bu_h(t) \|_{L^1} \le | \Omega |^{{\frac12}} \| \bu(t) - \bu_h(t) \|_{L^2}$, we have that for any $t$,
\begin{equation}\label{mom_lower}
\| \bu(t) - \bu_h(t) \|_{L^2} \ge  \max_i \frac{e_i^{mom}(t)}{|\Omega |^{{\frac12}} }
\end{equation}
Hence we observe that the momentum error serves as a lower bound for the $L^2$ velocity error at any time. By the same argument, the deviation in the angular momentum provides a similar lower bound  on the velocity error, but with a different constant in the denominator.

Of course, incorrect energy prediction also prevents convergence of $\bu_h$ to $\bu$, which is revealed by the estimate
\begin{align*}
E(\bu)-E(\bu_h)& = \frac12 \| \bu\|^2 - \|\bu_h \|^2
 = \frac12(\|\bu\|_{L^2}-\|\bu_h\|_{L^2})(\|\bu\|_{L^2}+\|\bu_h\|_{L^2})\\
& \le
\frac12(\|\bu-\bu_h\|_{L^2})(\|\bu\|_{L^2}+\|\bu_h\|_{L^2})
\le \frac1{\sqrt{2}}(\|\bu-\bu_h\|_{L^2})(E(\bu)+E(\bu_h))^{\frac12}.
\end{align*}
This implies that
\begin{equation}\label{energy_lower}
\| \bu(t) - \bu_h(t) \|_{L^2} \ge  \frac{\sqrt{2}|e^{E}(t)|}{(E(\bu)+E(\bu_h))^{\frac12}},
\end{equation}
where $e_E(t) = (E(\bu)-E(\bu_h))(t)$ is the energy error at time $t$.

This is evidence that the EMAC formulation, whose
nonlinear term correctly does not contribute to the momentum, angular momentum and energy balances, would be expected to (in general) not have as large of a lower bound on its $L^2$ velocity error as would a scheme whose nonlinearity nonphysically and incorrectly contributes to the balances.

\subsection{An error estimate with $Re$-independent Gronwall constant}

We now compare convergence estimates of EMAC and SKEW.  The result for SKEW is now considered classical, and can be found in, e.g. \cite{laytonbook,temam}. For comparison purposes, and since the EMAC scheme analysis is the same as SKEW except for the nonlinear terms, we first give the result for SKEW along with a brief proof following the exposition in \cite{J16}, which we found to be the most straightforward and concise for our purposes herein.  Afterwards, we prove a result for EMAC.  We will see that the fundamental difference is in the constants arising in the two schemes after the Gronwall inequality is applied.  We note that a convergence result for ROT will be the same as SKEW except with usual pressure is replaced by Bernoulli pressure, and no such results are known for CONV and CONS, as they are not energy stable (due to lack of energy conservation) and the arguments used in the proof below will no longer hold.


\begin{Theorem}\label{skewconv} [Convergence of SKEW]
Let $\bu_h$ solve \eqref{skew1}--\eqref{skew2} and $(\bu,p)$ be an NSE solution with $\bu_t\in L^2(0,T;\bX')$, $\bu\in L^4(0,T;\bH^1(\Omega))$, and $p\in L^2(0,T;L^2(\Omega))$.  Denote $\be(t)=\bu(t)-\bu_h(t)$ and $\bfeta(t)=\bu(t)-I_{St}^h \bu(t)$.\\ (i) It holds for all $t$ in $(0,T]$,
\begin{multline}
\| \be(t) \|^2  + \nu \int_0^t\| \nabla \be \|^2dt \le C \bigg( \| \bfeta(t) \|^2 + \nu \| \nabla \bfeta \|_{L^2(0,t;\bL^2)}^2  +K \bigg(
\nu^{-1} \| \bfeta_t \|^2_{L^2(0,T;\bX')}\\  +\nu^{-1} \inf_{q_h\in L^2(0,t;Q_h)} \| p-q_h\|^2_{L^2(0,t;L^2)}
+{\color{black} C_s(\bu,\blf,\nu^{-1})}\| {\bfeta} \|_{L^4(0,t;\bH^1)}^2 \bigg) \bigg).
\label{EstScew}
\end{multline}
with
\begin{equation}\label{K1}
K= \exp\left( C\nu^{-3} \| \bu \|_{L^4(0,t;\bH^1)} \right)
\end{equation}
and {\color{black}  $C_s(\bu,\blf,\nu^{-1})$ is} a factor depending (polynomially) on $\nu^{-1}$, $\| \bu \|_{L^4(0,t;\bH^1)}$, $\| \bu_0 \|$ and $\| \blf \|_{L^2(0,t;\bX')}$.\\
(ii) If additionally $\nabla \bu\in L^1(0,T;\bL^{\infty}(\Omega))$ and $\bu\in L^2(0,T;\bL^{\infty}(\Omega))$, then \eqref{EstScew} holds with
\begin{equation}\label{K2}
K=\exp\left( C( \| \bu \|_{L^{{2}}(0,t;\bL^{\infty})} + \nu^{-1} \| \nabla \bu \|_{L^{{1}}(0,t;\bL^{\infty})}) \right).
\end{equation}
\end{Theorem}
The factor $K$, which scales many of the right hand side terms in the error bound, is critical.  Since $\nu$ is inversely proportional to the Reynolds number, when Reynolds numbers are large, or even moderate, $K$ can be sufficiently large to render the error bound useless.  The term arises from the application of the Gronwall inequality applied to a term resulting from the analysis of the skew-symmetric nonlinearity in \eqref{ss8}.  We give a sketch of the proof below, and refer interested readers to \cite{J16} if more details are desired.

\begin{proof}

Setting $\be=\bu-\bu_h$ and using \eqref{nse1}--\eqref{nse2}, \eqref{skew1}--\eqref{skew2}, an error equation arises of the form
\begin{equation}
(\be_t,\bv_h) + b^*(\bu,\bu,\bv_h) - b^*(\bu_h,\bu_h,\bv_h) - (p,\nabla \cdot \bv_h) + \nu (\nabla \be,\nabla \bv_h)  = 0 \quad \forall \bv_h \in \bV_h. \label{ss2}
\end{equation}
Decomposing the error as $\be=(\bu-I_{St}^h \bu) + (I_{St}^h \bu - \bu_h) =: \bfeta + \bphi_h$, we obtain from \eqref{ss2} that for any $q_h\in Q_h$ and $\bv_h\in\bV_h$, it holds
\begin{multline}\label{ss3}
(\bfeta_t,\bv_h) + (({\bphi_h})_t,\bv_h) + b^*(\bu,\bu,\bv_h) - b^*(\bu_h,\bu_h,\bv_h)  \\
- (p - q_h,\nabla \cdot \bv_h) + \nu (\nabla \bfeta,\nabla \bv_h) +  \nu (\nabla \bphi_h,\nabla \bv_h) = 0.
\end{multline}
Now choosing $\bv_h=\bphi_h$ vanishes 
the last viscous term thanks to definition of the discrete Stokes projection, and yields
\begin{align}
\frac12 \frac{d}{dt} & \| \bphi_h \|^2 + \nu \| \nabla \bphi_h \|^2 =  {\color{black} -  (\bfeta_t,\bphi_h)}
-b^*(\bu,\bu,\bphi_h) + b^*(\bu_h,\bu_h,\bphi_h)- (p - q_h,\nabla \cdot \bphi_h). \label{ss4}
\end{align}
It is convenient to split the remainder of the proofs in three steps.\\
Step 1.
The linear terms are majorized using the $\bX'$ norm, Cauchy-Schwarz and Young inequalities:
\begin{align*}
\bigg|  (\bfeta_t,\bphi_h) \bigg| & \le 2\nu^{-1} \| \bfeta_t \|_{\bX'}^2 + \frac{\nu}{8} \| \nabla \bphi_h \|^2, \\
\bigg| (p - q_h,\nabla \cdot \bphi_h) \bigg| & \le  2\nu^{-1} \| p-q_h \|^2 + \frac{\nu}{8} \| \nabla \bphi_h\|^2.
\end{align*}
To bound the nonlinear terms, we split the difference of two trilinear forms as follows:
\begin{equation}\label{aux348}
b^*(\bu,\bu,\bphi_h) - b^*(\bu_h,\bu_h,\bphi_h)=b^*(\bfeta,\bu,\bphi_h)+ b^*(\bphi_h,\bu,\bphi_h) + b^*(\bu_h,\bfeta,\bphi_h).
\end{equation}
Application of H\"older's inequality, Sobolev embeddings, and Young's inequality produce the bounds:
\begin{align}
 \bigg| b^*(\bfeta,\bu,\bphi_h) \bigg| & \le C\nu^{-1} \| \nabla \bu \|^2 \|  \bfeta \| \| \nabla \bfeta \| +  \frac{\nu}{8} \| \nabla \bphi_h \|^2, \label{ss6} \\
  \bigg| b^*(\bu_h,\bfeta,\bphi_h) \bigg| & \le C\nu^{-1} \| \bu_h \| \| \nabla \bu_h \| \| \nabla \bfeta \|^2 +  \frac{\nu}{8} \| \nabla \bphi_h \|^2, \label{ss7} \\
  \bigg| b^*(\bphi_h,\bu,\bphi_h)  \bigg| & \le \frac12 \| \bphi_h \|_{L^3} \| \nabla \bu\| \| \bphi_h \|_{L^6} +  \frac12 \| \bphi_h \|_{L^3} \| \nabla \bphi_h\| \| \bu \|_{L^6}  \nonumber \\
  & \le C \| \nabla \bu \| \| \bphi_h \|^{1/2} \| \nabla \bphi \|^{3/2} \nonumber \\
& \le C \nu^{-3} \| \nabla \bu \|^4 \| \bphi_h \|^2 + \frac{\nu}{8} \| \nabla \bphi_h \|^2, \label{ss8}
\end{align}
with $C$ being a generic constant depending only on the domain.

Step 2. Combining the bounds now provides
\begin{align}
 \frac{d}{dt} & \| \bphi_h \|^2 + \frac{3\nu}{8} \| \nabla \bphi_h \|^2 \le
  2\nu^{-1} \| \bfeta_t \|_{\bX'}^2
 + 2\nu^{-1} \| p-q_h \|^2
 \nonumber \\
 &
+  C\nu^{-1} \left( \| \nabla \bu \|^2 \| \bfeta \| \| \nabla \bfeta \| + \| \bu_h \| \| \nabla \bu_h \| \| \nabla \bfeta \|^2 \right)
 + C \nu^{-3} \| \nabla \bu \|^4  \| \bphi_h \|^2
. \label{ss9}
\end{align}
Now integrating \eqref{ss2} over $(0,t)$, appropriately bounding the first two nonlinear terms using Cauchy-Schwarz and the usual stability estimates (e.g. from \cite{J16,laytonbook}), applying the Gronwall inequality, and using that $\bphi_h(0)={\bf 0}$ produces for any $t$ in $(0,T]$,
\begin{multline}
\|  \bphi_h(t) \|^2  + \frac{3\nu}{8} \| \bphi_h \|_{L^2(0,t;\bH^1)}^2  \le
C\nu^{-1} \exp \left( C \nu^{-3} \| \bu \|_{L^4(0,t;\bH^1)} \right) \bigg(
 \| \bfeta_t \|_{L^2(0,t;\bX')}^2\\  +\| p-q_h\|^2_{L^2(0,t;L^2)}
{\color{black} +C_s(\bu,\blf,\nu^{-1})}\| \bfeta \|_{L^4(0,t;\bH^1)} \bigg). \label{ss10}
\end{multline}
From here,  applying the triangle inequality completes the proof of \eqref{EstScew} with $K$ from \eqref{K1}.\\

Step 3. To show \eqref{EstScew} with $K$ from \eqref{K2}, one only has to repeat steps 1 and 2  with the only changes made to \eqref{ss8}:  Apply  integration by parts to the second term of the skew-symmetric form, H\"older's inequality, and Young's inequality to obtain
\begin{multline} \label{ss21}
\bigg| b^*(\bphi_h,\bu,\bphi_h) \bigg| = \bigg| \frac12 b(\bphi_h,\bu,\bphi_h) - \frac12 b(\bphi_h,\bphi_h,\bu) \bigg|\\ = \bigg|  b(\bphi_h,\bu,\bphi_h)  + \frac12 ((\nabla \cdot \bphi_h)\bphi_h,\bu) \bigg|
 \le \| \nabla \bu \|_{\bL^{\infty}} \| \bphi_h \|^2 + \frac12 \| \bu \|_{\bL^{\infty}} \| \nabla \cdot \bphi_h\| \| \bphi_h \|\\
   \le (\| \nabla \bu \|_{\bL^{\infty}}  +   \nu^{-1} \| \bu \|_{\bL^{\infty}}^2 )    \| \bphi_h \|^2 + \frac{\nu}{8} \| \nabla \bphi_h \|^2.
\end{multline}
\end{proof}

\begin{remark}\rm
An improvement offered by \eqref{K2} compared to \eqref{K1}  is that the Gronwall constant exponent has a $\nu^{-1}$ instead of a $\nu^{-3}$ for the price of higher regularity assumption.  It does not appear possible to completely remove an inverse dependence on $\nu$ from the Gronwall constant.  From \eqref{ss21} one soon notes that the explicit $\nu$-dependence of the exponent factor would disappear if the finite element solution satisfies $\nabla\cdot\bu_h=0$ pointwise {\cite{SLLL18}, or if grad-div stabilization is utilized \cite{FGJN18}} . This provoked an opinion that a stronger enforcement of div-free is necessary for ``robust'' error estimates in mixed FE methods for the Navier--Stokes equations. The theorem coming next shows that the EMAC formulation removes such dependence for  standard mixed elements, like Taylor-Hood, {without any divergence stabilization}, and delivers the same  Gronwall factors in \eqref{EstScew} as pointwise divergence-free elements in~\cite{SLLL18}. This supports our hypothesis that the preservation of energy and kinematic balances is crucial, while using div-free elements is a possible (although sometimes expensive) way of ensuring these balances for FE formulation with other forms of nonlinear terms.
\end{remark}

We consider now the analogous convergence result for the EMAC scheme.

\begin{Theorem}\label{emacconv} [Convergence of EMAC]
Let $\bu_h$ solve \eqref{emac1}--\eqref{emac2} and $(\bu,p)$ be an NSE solution with $\bu_t\in L^2(0,T;\bX')$, $\bu\in L^4(0,T;\bH^1(\Omega))$,  $\nabla\bu\in L^2(0,T;\bL^3(\Omega))$, and $P=(p-\frac12 |\bu|^2)\in L^2(0,T;L^2(\Omega))$.
Denote $\be(t)=\bu(t)-\bu_h(t)$ and $\bfeta(t)=\bu(t)-I_{St}^h \bu(t)$. \\
(i) For all $t$ in $(0,T]$, it holds
\begin{multline}
\| \be(t) \|^2  + \nu \int_0^t\| \nabla \be \|^2dt \le C \bigg( \| \bfeta(t) \|^2 + \nu \| \nabla \bfeta \|_{L^2(0,t;\bL^2)}^2  +K \bigg(
\nu^{-1} \| \bfeta_t \|^2_{L^2(0,T;\bX')}\\  +\nu^{-1} \inf_{q_h\in L^2(0,t;Q_h)} \| P-q_h\|^2_{L^2(0,t;L^2)}
+{\color{black} C_e}{(\bu,\nu^{-1})} \| \bfeta \|_{L^4(0,t;\bH^1) }\bigg) \bigg).
\label{EstEmac}
\end{multline}
with
\begin{equation}\label{K1emac}
K= \exp\left( C\nu^{-1} \| \nabla\bu \|_{L^2(0,t;L^3)}\right),
\end{equation}
and {\color{black} $C_e(\bu,\nu^{-1})$ is} a factor depending (polynomially) on $\nu^{-1}$ and $\| \bu \|_{L^4(0,t;\bH^1)}$.\\
(ii) If additionally $\bu, \nabla \bu\in L^1(0,T;\bL^{\infty}(\Omega))$, then \eqref{EstEmac} holds with
\begin{equation}\label{K2emac}
K=\exp\left( C( \| \bu \|_{L^1(0,t;\bL^{\infty})} +  \| \nabla \bu \|_{L^1(0,t;\bL^{\infty})} ) \right),
\end{equation}
and also with $C(\bu,\nu^{-1})=C (\| \bu \|_{L^1(0,t;\bL^{\infty})} +  \| \nabla \bu \|_{L^1(0,t;\bL^{\infty})})$, independent of $\nu^{-1}$ .
\end{Theorem}

\begin{remark}\rm
The key difference of the EMAC convergence  theorem compared to the result for SKEW is that the Gronwall constant $K$ from EMAC has no explicit dependence on the
viscosity, if the solution $\bu$ satisfies the same\footnote{Actually, assumptions here are even slightly weaker for EMAC.} additional regularity assumptions as in part (ii) of Theorem~\ref{skewconv}. Moreover, a factor multiplying $\| \bfeta \|_{L^4(0,t;\bH^1)}$ is also independent of $\nu$. It is interesting to see that for the basic regularity case in part (i) of the theorem, the $\nu$-dependence in the $K$ also reduces from $\nu^{-3}$ for SKEW to $\nu^{-1}$ for EMAC.
\end{remark}

{\color{black}
\begin{remark} \rm
 The dependence of velocity error on  pressure approximation in \eqref{EstEmac} (not present for div-free elements) can usually be largely ameliorated by the simple grad-div stabilization~\cite{O02,reusken}, or by element choices such as the mini element~\cite{arnold1984stable,GR86}, although for very large or complex pressures it may be necessary to use large grad-div stabilization parameters or divergence-free elements \cite{LR19} to sufficiently reduce numerical error arising from the pressure approximation.
 \end{remark}

 \begin{remark}\label{emacp} \rm
  \MO{Related to the above remark on pressure and velocity errors interference,  we note that} one can construct examples such as \MO{Poiseuille}  flow, where EMAC \MO{performs worse} than SKEW due to the EMAC pressure $P=p-\frac12|\bu|^2$ being larger \MO{or} more complex \MO{compared to the kinematic} pressure $p$. \MO{Nevertheless,} in most simulations of practical problems \MO{using EMAC} this does not happen, e.g. channel flow past a step and the Gresho problem in \cite{CHOR17}, a large variety of simulations done by Lehmkuhl and coworkers e.g. \cite{PCLRH18,LHOCR19,MSLGD19,LPH19} to name a few, and numerical tests below.  This is surprising since the Bernoulli pressure $p+\frac12 |\bu|^2$ is believed responsible for the generally poor accuracy of \MO{finite element codes using rotational formulation}.  While this \MO{phenomenon} is not completely understood, one possible explanation is that \MO{in (statistically) developed flows} pressure correlates with the square of the velocity with $\frac{p'}{\bu^{'^{2}}} \sim 1$ for low Reynolds numbers and $\frac{p'}{\bu^{'^{2}}} \sim 0.7$  \MO{for higher Reynolds numbers}~\cite{CC99}, where $'$ represents root mean square of the {fluctuation of the} quantity, and in this \MO{sense  $p-\frac12|\bu|^2$} may actually have reduced complexity compared to $p$.  This general lack-of-ill-effect of the EMAC pressure will be explored in future work of the authors.
 \end{remark}
}

\begin{proof}
The proof follows similar to that of Theorem \ref{skewconv} except for the trilinear forms, which now are $ c(\bu,\bu,\bv)=2(D(\bu)\bu,\bv)+(\mbox{div}(\bu)\bu,\bv)$ instead of $b^\ast(\bu,\bu,\bv)$. All arguments from the proof of Theorem~\ref{skewconv} before Step 1 remain literally the same, while in Step~1 in place of  \eqref{aux348} we use a slightly different decomposition for the difference  $c(\bu,\bu,\bphi_h)-c(\bu_h,\bu_h,\bphi_h)$. Namely, we decompose:
\begin{equation}\label{aux1}
\begin{aligned}
(\bD(\bu)\bu,\bphi_h)- &(\bD(\bu_h)\bu_h,\bphi_h)=(\bD(\bfeta)\bu,\bphi_h) +(\bD(I_{St}^h(\bu))\bu,\bphi_h) - \bD(\bu_h)\bu_h,\bphi_h)  \\
   =& (\bD(\bfeta)\bu,\bphi_h) + (\bD(I_{St}^h(\bu))\bfeta,\bphi_h)\\
    &-(\bD(I_{St}^h(\bu))\bphi_h,\bphi_h)-(\bD(\bphi_h) I_{St}^h(\bu),\bphi_h)- (\bD(\bphi_h)\bphi_h,\bphi_h),
\end{aligned}
\end{equation}
and similarly,
\begin{equation}\label{aux2}
\begin{aligned}
(\mbox{div}(\bu)\bu,\bphi_h)- &(\mbox{div}(\bu_h)\bu_h,\bphi_h)=(\mbox{div}(\bfeta)\bu,\bphi_h) + (\mbox{div}(I_{St}^h(\bu))\bfeta,\bphi_h)\\
& -(\mbox{div}(I_{St}^h(\bu))\bphi_h,\bphi_h)-(\mbox{div}(\bphi_h) I_{St}^h(\bu),\bphi_h)- (\mbox{div}(\bphi_h)\bphi_h,\bphi_h).
\end{aligned}
\end{equation}
Before combining \eqref{aux1} and \eqref{aux2}, we first manipulate the $(\bD(\bphi_h) I_{St}^h(\bu),\bphi_h)$ term in \eqref{aux1}. Using the definition $2\bD(\bphi_h)=\nabla\bphi_h+(\nabla\bphi_h)^T$ and integration by parts, we get
\begin{equation}\label{aux3}
\begin{aligned}
2(D(\bphi_h) &I_{St}^h(\bu),\bphi_h)= ((\nabla\bphi_h)I_{St}^h(\bu),\bphi_h)+ ((\nabla\bphi_h)\bphi_h,I_{St}^h(\bu))\\
&= (I_{St}^h(\bu)\cdot\nabla\bphi_h ,\bphi_h)+ (\bphi_h\cdot\nabla\bphi_h,I_{St}^h(\bu))\\
&=-\frac12(\mbox{div}(I_{St}^h(\bu))\bphi_h ,\bphi_h)- (\bphi_h\cdot\nabla I_{St}^h(\bu),\bphi_h) - (\mbox{div}(\bphi_h) I_{St}^h(\bu), \bphi_h)\\
&=-\frac12(\mbox{div}(I_{St}^h(\bu))\bphi_h ,\bphi_h)- (\bD(I_{St}^h(\bu))\bphi_h,\bphi_h) - (\mbox{div}(\bphi_h) I_{St}^h(\bu), \bphi_h).
\end{aligned}
\end{equation}

Summing up the equality \eqref{aux1} scaled by 2 with \eqref{aux2} and using \eqref{aux3} along with \eqref{skew} (noting that the term $- 2(\bD(\bphi_h)\bphi_h,\bphi_h)- (\mbox{div}(\bphi_h)\bphi_h,\bphi_h)$ vanishes thanks to EMAC's nonlinearity preserving energy), leads to
\begin{align*}
c(\bu,\bu,\bphi_h)-c(\bu_h,\bu_h,\bphi_h)&=
   ([2\bD(\bfeta)+\mbox{div}(\bfeta)]\bu,\bphi_h) + ([2\bD(I_{St}^h(\bu))+\mbox{div}(I_{St}^h(\bu))]\bfeta,\bphi_h)\\
   &- ([\bD(I_{St}^h(\bu))+\frac12\mbox{div}(I_{St}^h(\bu))]\bphi_h,\bphi_h).
\end{align*}
Now we estimate the terms on the right-hand side, depending on our regularity assumption for $\bu$. For part (i) of the theorem, we
repeat the estimates from \eqref{ss6} and \eqref{ss7} and use the  $H^1$-stability of the Stokes interpolant to get
\begin{multline}\label{ss6-7}
  |([2\bD(\bfeta)+\mbox{div}(\bfeta)]\bu,\bphi_h) + ([2\bD(I_{St}^h(\bu))+\mbox{div}(I_{St}^h(\bu))]\bfeta,\bphi_h)|\\
  \le C\nu^{-1}( \| \nabla \bu \|^2 \|  \bfeta \| \| \nabla \bfeta \| + \| \nabla \bu \| \|  \bu\| \| \nabla \bfeta \|^2) +  \frac{\nu}{4} \| \nabla \bphi_h \|^2
\end{multline}
and
\begin{multline}\label{ss8emac}
 |([\bD(I_{St}^h(\bu))+\frac12\mbox{div}(I_{St}^h(\bu))]\bphi_h,\bphi_h)|\le \frac32\|\bD(I_{St}^h(\bu)\|_{L^3}\|\bphi_h\|_{L^2}\|\bphi_h\|_{L^6}\\ \le
 C\|\nabla\bu\|_{L^3}\|\bphi_h\|_{L^2}\|\nabla\bphi_h\|\le  C\nu^{-1}\|\nabla\bu\|_{L^3}^2\|\bphi_h\|_{L^2}^2 +\frac\nu8\|\nabla\bphi_h\|^2.
\end{multline}
The second estimate above uses the stability of the Stokes projection in \eqref{stokesbound} with $r=3$. We now proceed as in step 2 of the proof of Theorem~\ref{skewconv} to show \eqref{EstEmac}--\eqref{K1emac}.

To prove part (ii) of the theorem, for  $\bu, \nabla \bu\in L^1(0,T;\bL^{\infty}(\Omega))$ we update estimates in \eqref{ss6-7}--\eqref{ss8emac} as follows:
\begin{align*}
|c(\bu,\bu,\bphi_h)-c(\bu_h,\bu_h,\bphi_h)|\le
   & \frac{9}{2}\|\bD(\bfeta)\|_{L^2}^2\|\bu\|_{L^\infty} + \frac{9}{2}\|\bD(I_{St}^h(\bu))\|_{L^\infty}\|\bfeta\|_{L^2}^2\\
   & +
   \left[\|\bu\|_{L^\infty} + \frac52\|\bD(I_{St}^h(\bu))\|_{L^\infty}\right]\|\bphi_h\|_{L^2}^2.
\end{align*}
Finally using  \eqref{stokesbound} with $r=\infty$ and reducing leads to
\begin{multline}\label{est}
  |c(\bu,\bu,\bphi_h)-c(\bu_h,\bu_h,\bphi_h)|\\ \le C\Big(\|\nabla \bfeta\|^2\|\bu\|_{\bL^\infty} + \|\nabla \bu \|_{\bL^\infty}\|\bfeta\|^2
    +
   \left[\|\bu\|_{\bL^\infty} + \|\nabla \bu \|_{\bL^\infty}\right]\|\bphi_h\|^2\Big).
\end{multline}
From here, the proof follows the arguments of step 2
in the proof of Theorem \ref{skewconv}.
\end{proof}

Error estimates \eqref{EstScew} and \eqref{EstEmac}  can be developed further in a standard way by using textbook (e.g., \cite{GR86})  approximation properties of the Stokes projection in Hilbert spaces.

{\color{black}

\section{Numerical Examples}

We present here results of numerical experiments that extend the numerical testing of EMAC.
All tests will use $(P_2,P_1)$ Taylor-Hood elements (with no stabilization) for velocity and pressure, and for temporal discretizations we use BDF2 in the schemes \eqref{emac1}-\eqref{emac2} for EMAC, and \eqref{skew1}-\eqref{skew2} for SKEW.  Although many comparisons of schemes and formulations could be done, we focus herein on this comparison, since SKEW is the most widely used scheme that preserves energy (which is widely believed critical for higher Reynolds number flows \cite{CHOR17,CHOR19}, and CONS and CONV do not preserve it); although ROT preserves energy, it suffers from numerical difficulties arising from use of Bernoulli pressure $P=(p+\frac12|\bu|^2)$ that must be stabilized either with divergence-free elements or (sometimes heavy) grad-div stabilization \cite{O02,LMNOR09,GRT14,JLMNR17,RKB17,OHRG15}, which makes it less attractive for general use; still, we include ROT in our first numerical test, and not surprisingly it performs poorly.  It is interesting, and in the first few EMAC papers was surprising, that the EMAC pressure  $P_{emac}=(p-\frac12|\bu|^2)$  does not have the ill-effect that Bernoulli pressure has on ROT, as can be seen in results below, see also Remark \ref{emacp}.

\subsection{Numerical test with known analytical solution}


The first test problem we consider is the lattice vortex problem from \cite{MB2002,SL17}.  The true solution takes the form
\[
\bu = \bv e^{-8\nu \pi^2 t}, \ \ p=q e^{-16 \pi^2 t},
\]
where
\[
\bv = \langle \sin (2\pi x) \sin(2\pi y),\ \cos(2\pi x)\cos(2\pi y) \rangle,\ \ \ q = -\frac12 \left( \sin^2(2 \pi x) + \cos^2(2 \pi y) \right).
\]
Notice that $(\bu,p)$ is an exact NSE solution with $\blf={\bf 0}$, and that the velocity has \MO{zero} momentum and angular momentum at all times.  This problem was
used in \cite{SLLL18} to illustrate how the improvement in the Gronwall constant in the error analysis for strongly divergence-free approximations leads to
more accurate solutions (than weakly divergence-free approximations) over longer time intervals, and was also used in \cite{CHOR19} to illustrate advantages of EMAC
over some other typical formulations.  Here, we use the test to compare SKEW, EMAC and ROT using $\nu=10^{-5}$ on [0,10].
For comparison, we also \MO{run the test for Scott-Vogelius (SV) elements and CONV formulation} (note that with SV elements, all of these formulations are equivalent).

We consider the \MO{lattice vortex} problem \MO{in} $\Omega=(0,1)^2$, with $\bu_0=\bu(0)$ as the initial condition and $\blf={\bf 0}$.
This is a difficult problem because there are several spinning vortices whose edges touch, which can be difficult to resolve numerically.
This problem is chosen since it has a known true solution and allows us to  calculate errors and make comparisons.
We simulate this problem using \MO{the} Crank-Nicolson time stepping \MO{scheme} with $\Delta t=0.001$, and $(P_2,P_1)$ elements on an Alfeld split of a Delaunay (non-uniform) mesh, with $h\sim \frac{1}{64}$.    This mesh \MO{fails to resolve all spatial scales} for this problem, and that is representative of the typical situation in Navier-Stokes simulations -- one can \MO{rarely have} a fine enough mesh in practice.  The nonlinear problems at each time step were solved with Newton's method, and Dirichlet boundary conditions were strongly enforced to be the true velocity solution at the boundary nodes.

\begin{figure}[!h]
\begin{center}
\includegraphics[width=.4\textwidth,height=0.25\textwidth, viewport=0 0 530 410, clip]{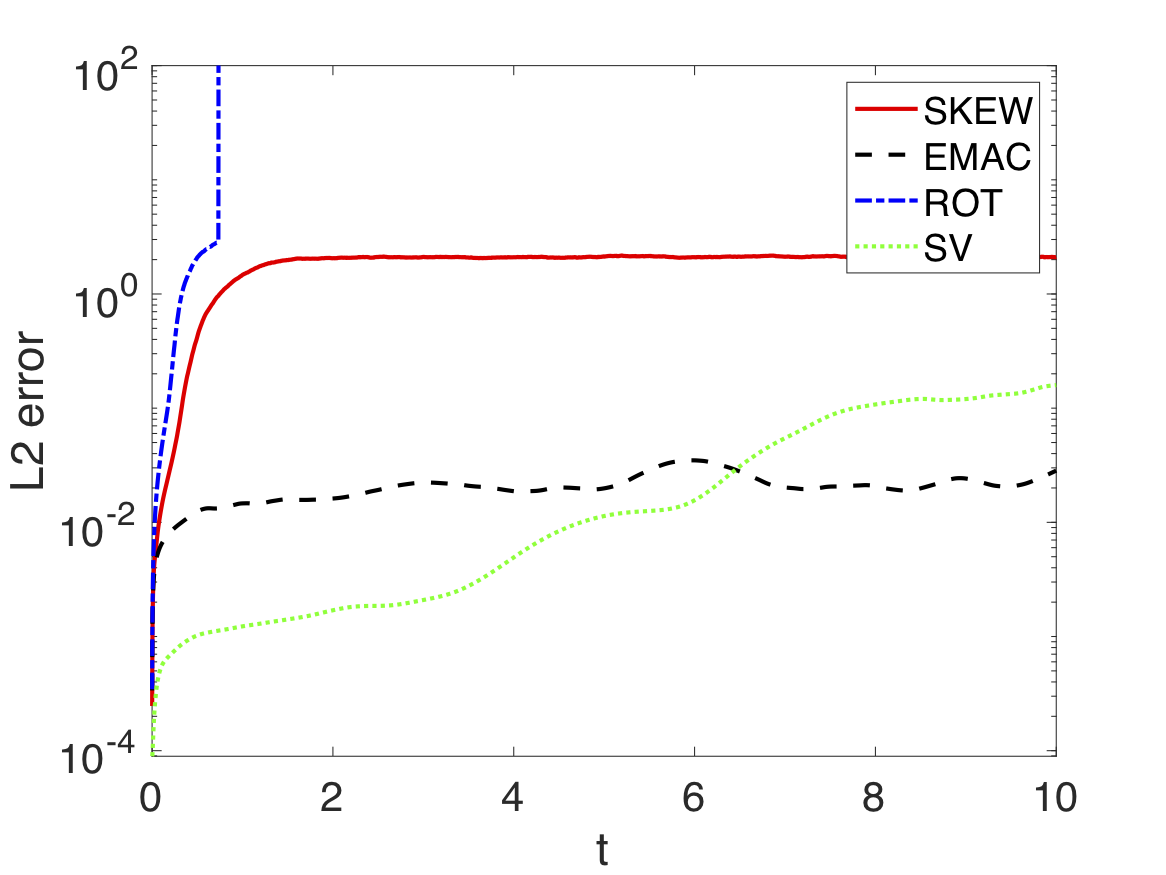}
\includegraphics[width=.4\textwidth,height=0.25\textwidth, viewport=0 0 530 410, clip]{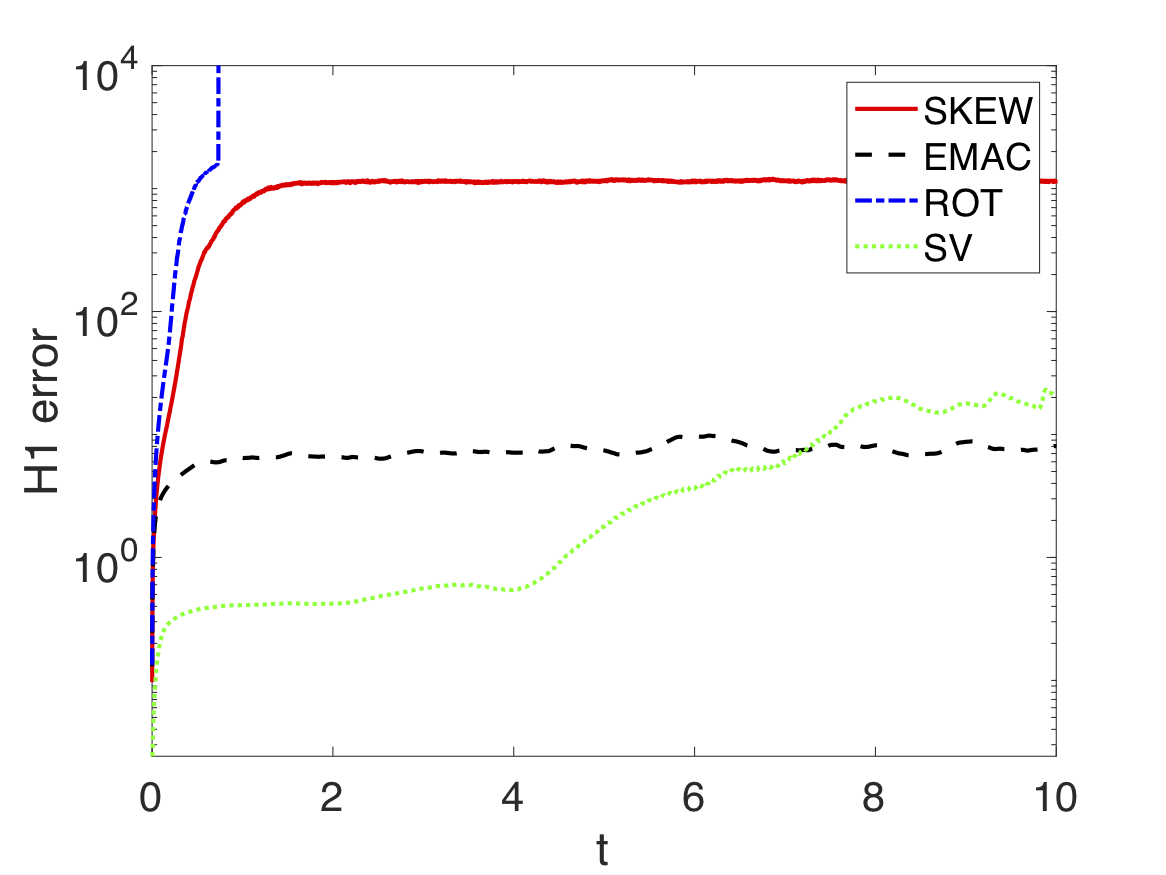}\\
\includegraphics[width=.4\textwidth,height=0.25\textwidth, viewport=0 0 530 410, clip]{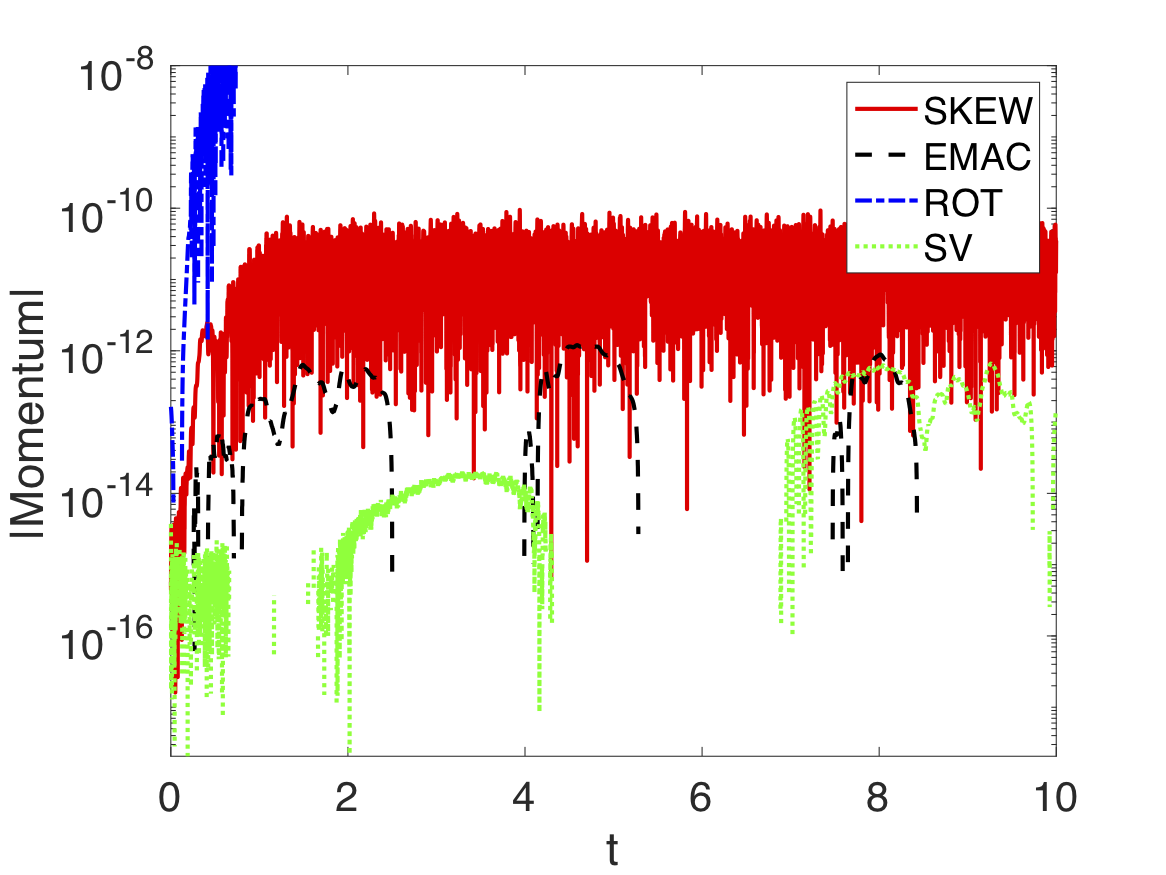}
\includegraphics[width=.4\textwidth,height=0.25\textwidth, viewport=0 0 530 410, clip]{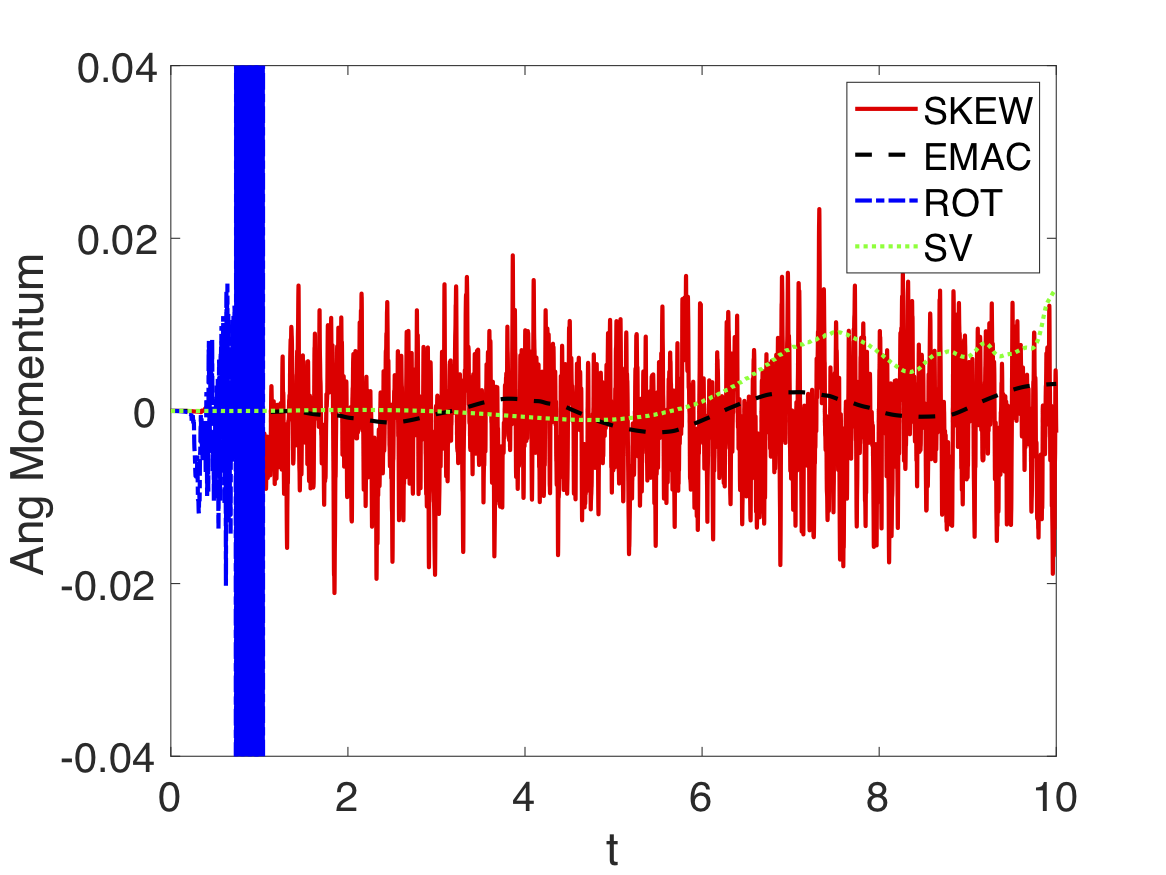}
\end{center}
\caption{\label{lattice1}
Shown above is the solution for the lattice vortex problem as a speed contour plot (left) and
pressure contour plot (right). }
\end{figure}

We show results of $L^2$ error, $H^1$ error, momentum, and angular momentum versus time, for each formulation in figure \ref{lattice1}.  Observe first that ROT, as expected, gives a very poor approximation, with error that is orders of magnitude worse than the other methods.  The advantage of EMAC over SKEW is clear: EMAC's error in $L^2$ and $H^1$ is 2 orders of magnitude smaller than for SKEW, and while both have momentum essentially at 0, the angular momentum of SKEW oscillates wildly while EMAC's is a smooth curve staying close to 0.  The SV solution is the most accurate at early times, but  for $t>7$, the SV error increases above that of EMAC and remains higher through $t=10$.  The SV solution's momentum and angular momentum are similar to those of EMAC, which is expected since the pointwise divergence-free property of SV (the $L^2$ divergence error remains below $10^{-10}$ at all times for SV), it will also preserve energy, momentum and angular momentum, just as EMAC does.

\subsection{2D channel flow past a cylinder}

For our next test we consider simulations for 2D channel flow past a cylinder, and compare results from EMAC and SKEW to each other and to a finer mesh solution.
This test problem was considered for comparing these formulations (and others) in \cite{CHOR17} for the case of time dependent inflow (max inlet velocity of $\sin(\pi t/8)$ on [0,8]) with $0\le Re(t) \le 100$, and essentially no difference was found between the formulations' solutions.  Here we consider the case of constant inflow with $Re=200$, which is a somewhat more challenging problem.

For the problem setup, we follow \cite{ST96,J04,MRXI17} which uses a $2.2\times0.41$ rectangular channel domain containing a cylinder (circle) of radius $0.05$ centered at $(0.2,0.2)$. There is no external forcing ($\blf={\bf 0}$), the kinematic viscosity is taken to be $\nu=0.0005$, no-slip boundary conditions are prescribed for the walls and the cylinder, while the inflow and outflow velocity profiles are given by
\begin{align*}
u_1(0,y,t) & = u_1(2.2,y,t) = \frac{6}{0.41^2}y(0.41-y),
\\u_2(0,y,t) & = u_2(2.2,y,t) = 0.
\end{align*}
All tests are run on $0\le t\le 10$, and start from rest, i.e. $u_0=0$.  A statistically steady, time-periodic state is reached by approximately $t=4$, and we calculated max and min lift and drag values over the time period $7\le t \le 10$.

\begin{figure}[!ht]
	\centering
		Mesh 1 \hspace{1.8in} Mesh 2\\
\includegraphics[width=.48\textwidth, height=.11\textwidth,viewport=140 20 1030 190, clip]{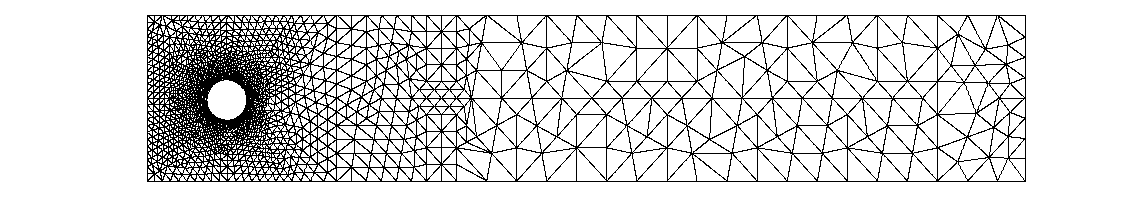}
	\includegraphics[width=.48\textwidth, height=.11\textwidth,viewport=140 20 1030 190, clip]{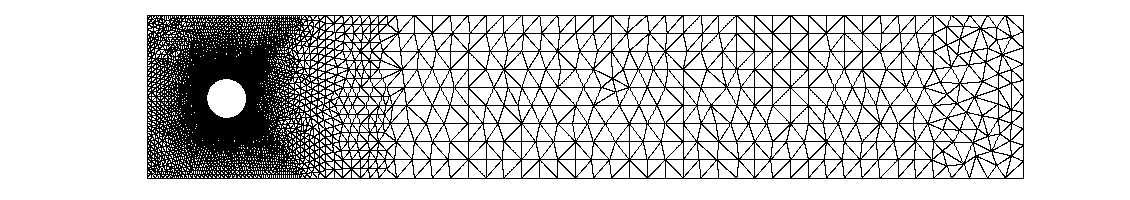}
	\caption{\label{meshes} Shown above are  the meshes used for the flow past a cylinder test problem.}
\end{figure}

Two meshes were used in our numerical tests, and are displayed in figure \ref{meshes}.  Mesh 2 is obtained by barycentric refinement of Mesh 1 and it is used only for the reference simulation; when equipped with $P_2$ velocity elements it provides 103K total velocity degrees of freedom (dof).  We compute on this mesh both with $(P_2,P_1^{disc})$ Scott-Vogelius (SV) elements and with SKEW using  $(P_2,P_1)$  Taylor-Hood (TH) elements and grad-div stabilization, using BDF2 time stepping and $\Delta t=0.001$.  Plots of these two fine mesh solutions are identical at all times, and their max lift and drag values agree to 4 digits.  Hence these fine mesh solutions \MO{are 
sufficiently} resolved to use as reference solutions for coarse mesh approximations; we will refer to the SV solution as the reference solution.

Mesh 1 is used to compare EMAC to SKEW, using $(P_2,P_1)$ Taylor-Hood elements, and it provides 35K total velocity dof.   Computations were run with BDF2 time stepping and $\Delta t=0.001$, and we note \MO{that} tests with $\Delta t=0.002$ gave very similar results.  Plots of solutions at $t=10$ are shown in figure \ref{Re200speed}, and we observe that all three plots appear resolved from the inlet up to the cylinder, which is not surprising since Mesh 1 is rather fine in this region.  However, in the downstream region where Mesh 1 is coarser, we observe EMAC provides a qualitatively accurate velocity prediction (although with some minor oscillations), while SKEW's \MO{solution is completely wrong} as it is destroyed by significant downstream oscillations that incorrectly predict a chaotic turbulent-like flow instead of the smoother vortex street.

\begin{figure}[!h]
\begin{center}
Reference  \\
\includegraphics[width=.45\textwidth, height=.15\textwidth,viewport=100 10 920 260, clip]{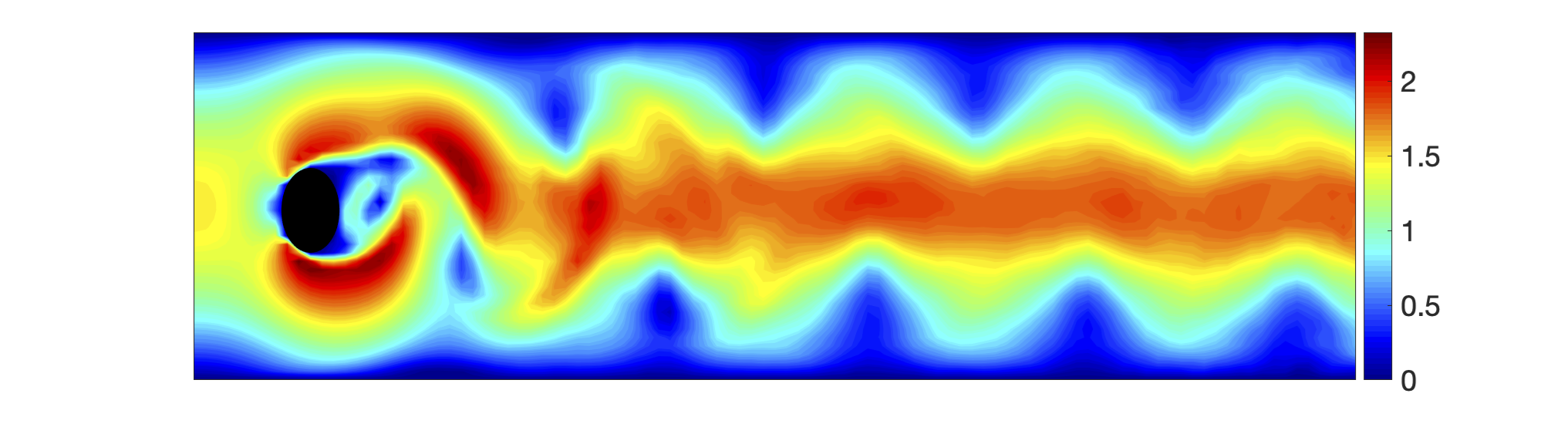}\\
EMAC   \hspace{1.75in} SKEW \\
\includegraphics[width=.45\textwidth, height=.15\textwidth,viewport=100 10 920 260, clip]{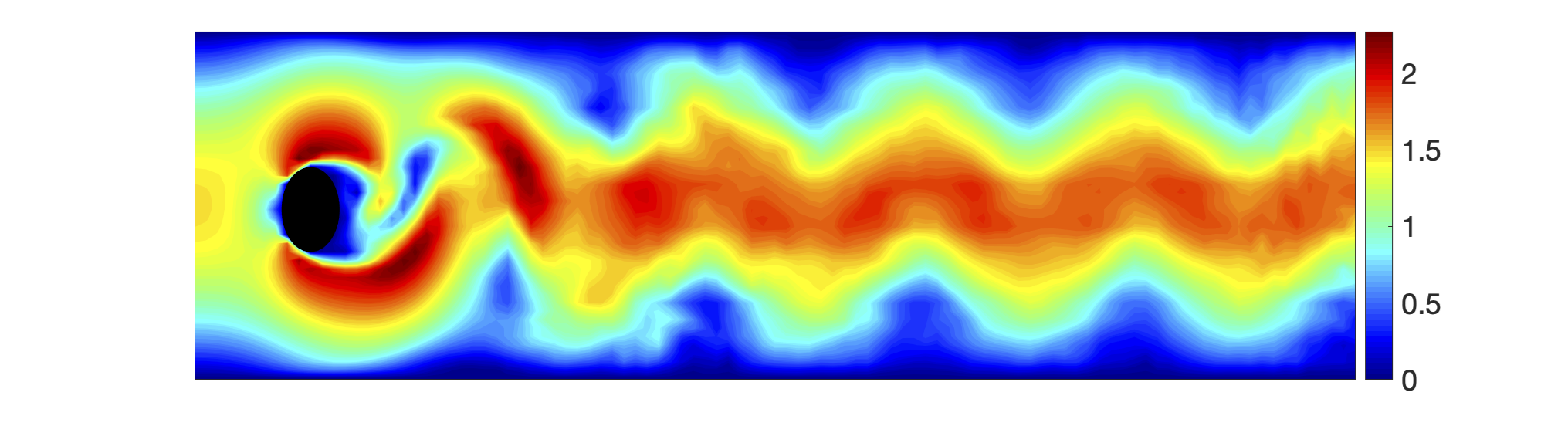}
\includegraphics[width=.45\textwidth, height=.15\textwidth,viewport=100 10  920 260, clip]{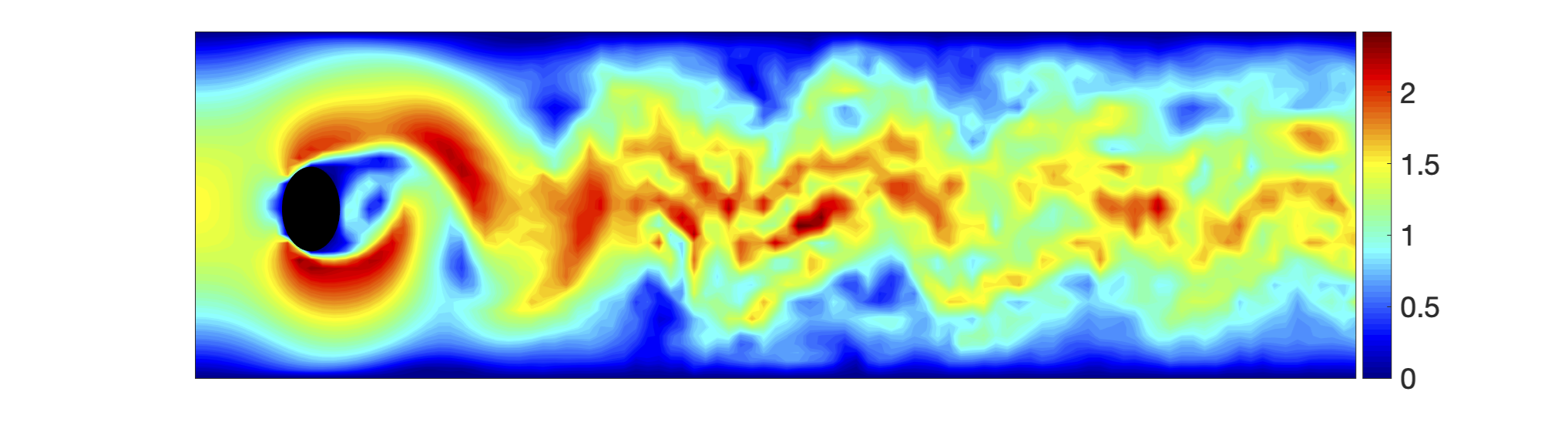}
\end{center}
\caption{\label{Re200speed}
Shown above are $Re=200$ solution speed contours at $t=10$ for DNS (SV on Mesh 2), EMAC (Mesh 1) and SKEW (Mesh 1).}
\end{figure}

We computed lift and drag values for the various solutions.  At each time step, we calculate lift and drag coefficients using global integral formulas derived following \cite{J04}, as they are less sensitive to the approximation of the circular boundary than boundary integrals are.
Let $\bv_d \in P_2(\tau_h)\cap C^0(\Omega)$ with $\bv_d |_{S} = (1,0)^T$ and that vanishes on all other boundaries, and similarly for $\bv_l$ but with $\bv_d |_{S} = (0,1)^T$.  There are multiple way to construct such functions $\bv_d,\bv_l$, and we chose to create them by a discrete Stokes extension of their respective boundary conditions.  The calculations are then done via:  for SKEW,
\begin{align*}
c_d(t^n) &= -20\left( \frac{1}{2\Delta t}(3\bu_h^{n} - 4\bu_h^{n-1} + \bu_h^{n-2},\bv_d) + b^*(\bu_h^{n},\bu_h^n,\bv_d) - (p^n,\nabla \cdot \bv_d) + \nu(\nabla \bu_h^n,\nabla \bv_d) \right), \\
c_l(t^n) &= -20\left(\frac{1}{2\Delta t}(3\bu_h^{n} - 4\bu_h^{n-1} + \bu_h^{n-2},\bv_l) + b^*(\bu_h^{n},\bu_h^n,\bv_l) - (p^n,\nabla \cdot \bv_l) + \nu(\nabla \bu_h^n,\nabla \bv_l) \right),
\end{align*}
and similar for EMAC with the SKEW nonlinearity replaced by the EMAC nonlinearity.  Results are shown in table \ref{data2DB}, and we observe EMAC and SKEW to have similar accuracy: EMAC has better {\color{black} lift} prediction, while SKEW has better {\color{black} drag} prediction.  The competitiveness of SKEW with EMAC for lift and drag calculations is due to Mesh 1 being sufficiently fine in the upstream region and around the cylinder -- with a sufficiently fine mesh, the formulations SKEW, EMAC, ROT, and CONV will all give accurate results.

\begin{table}[h!]
\centering
\begin{tabular}{|l|c|c|c|c|}
\hline
Method & $c_l^{max}$ &  $c_l^{min}$ 	&  $c_d^{max}$ &  $c_d^{min}$  \\ \hline
Mesh 2 (TH, SKEW + graddiv)      & 2.14511 & -2.19523 & 3.29090 & 2.96432 \\ \hline
Mesh 2 (SV)      & 2.14404 & -2.19422 & 3.29116 & 2.97689\\ \hline \hline  
Mesh 1 EMAC	& 2.12450	& -2.16637 &  3.30708 & 2.98819 \\ \hline
$|$ EMAC - SV $|$ & 1.95e-2 & 2.79e-2 & 1.59e-2 & 1.11e-2 \\ \hline \hline
Mesh 1 SKEW  & 2.11124 & -2.16652 & 3.28318 & 2.97413 \\ \hline
$|$ SKEW - SV $|$  & 3.28e-2 & 2.77e-2 & 7.98e-3 & 2.76e-3 \\ \hline
\end{tabular}
\caption{ \label{data2DB} \MO{Max} and min values for lift and drag, for the Reference, EMAC and SKEW \MO{solutions}, along with the errors of EMAC and SKEW compared to \MO{the} Reference.}
\end{table}

\subsection{2D Kelvin-Helmholtz simulation}

\begin{figure}[!h]
\begin{center}
EMAC \hspace{1.6in} SKEW \\
\includegraphics[width=.4\textwidth, height=.18\textwidth,viewport=0 0 520 395, clip]{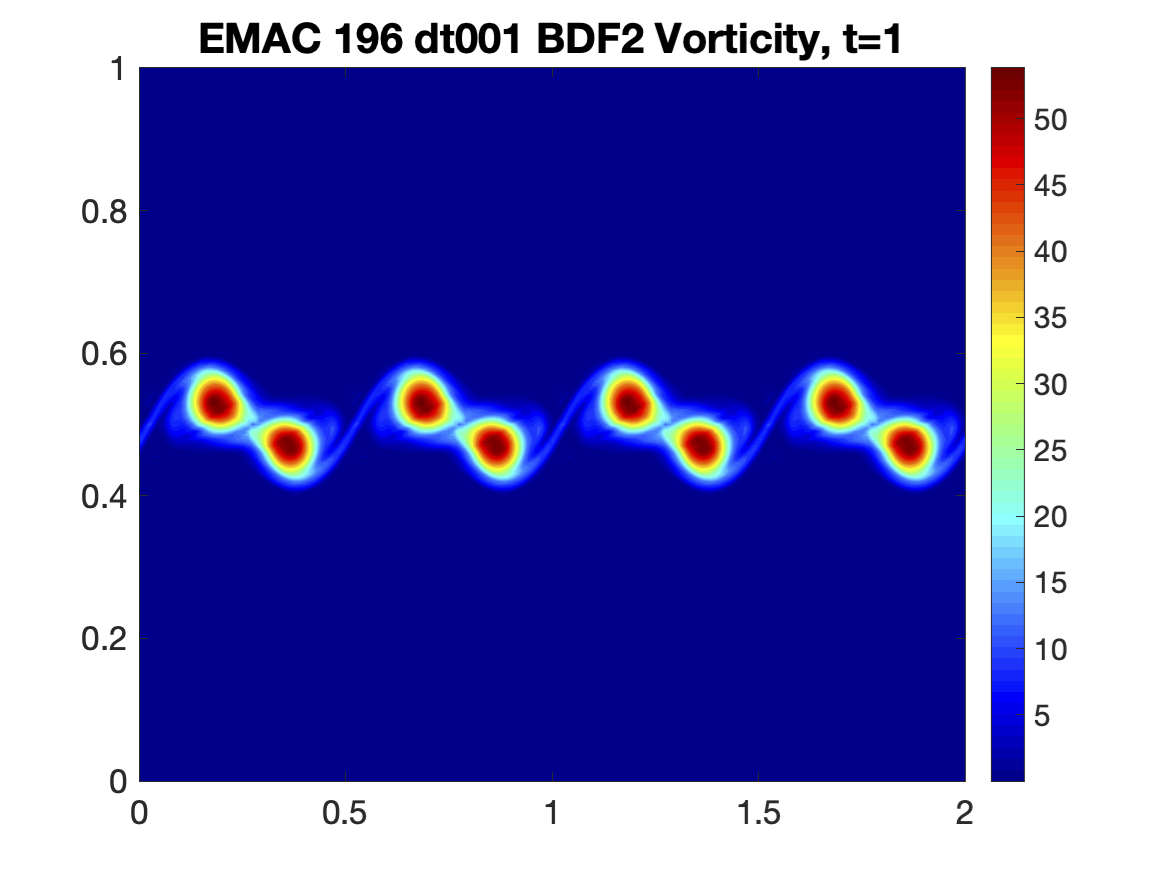}
\includegraphics[width=.4\textwidth, height=.18\textwidth,viewport=0 0 520 395, clip]{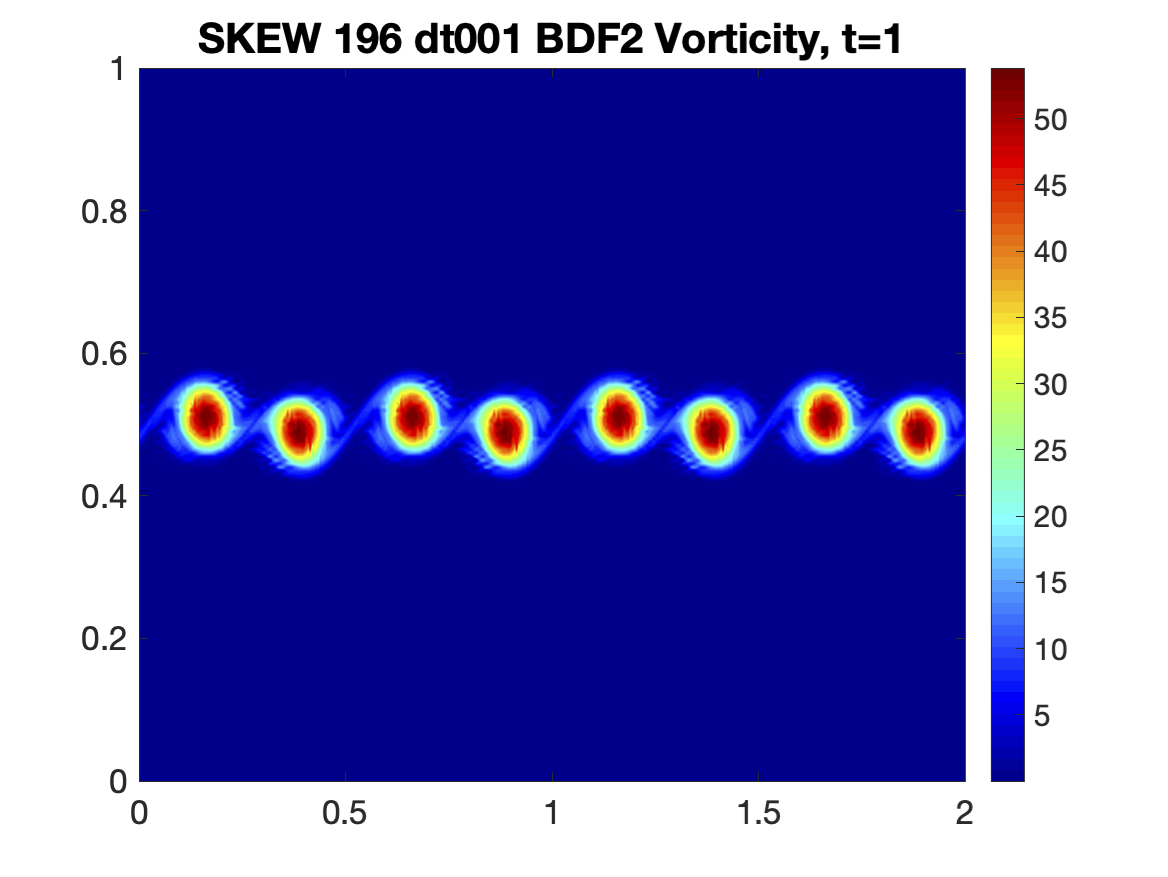} \\
\includegraphics[width=.4\textwidth, height=.18\textwidth,viewport=0 0 520 395, clip]{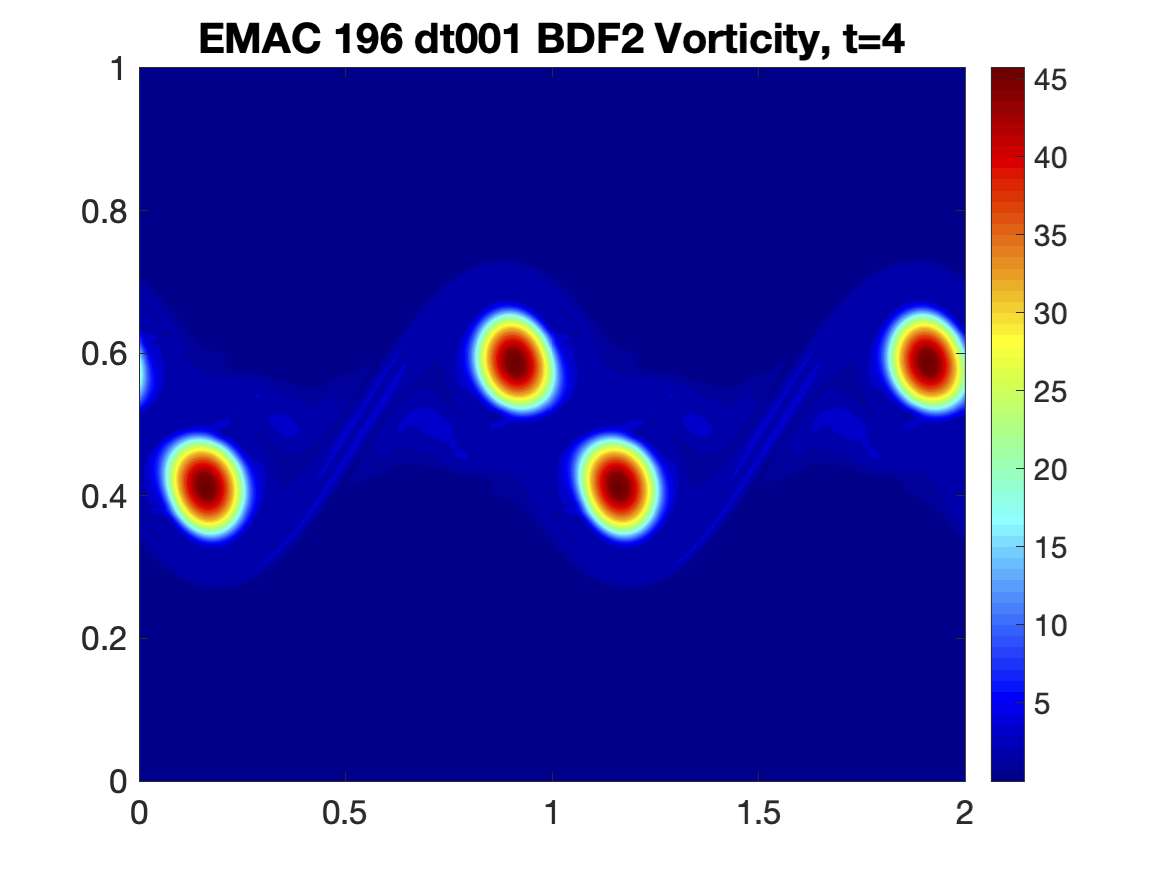}
\includegraphics[width=.4\textwidth, height=.18\textwidth,viewport=0 0 520 395, clip]{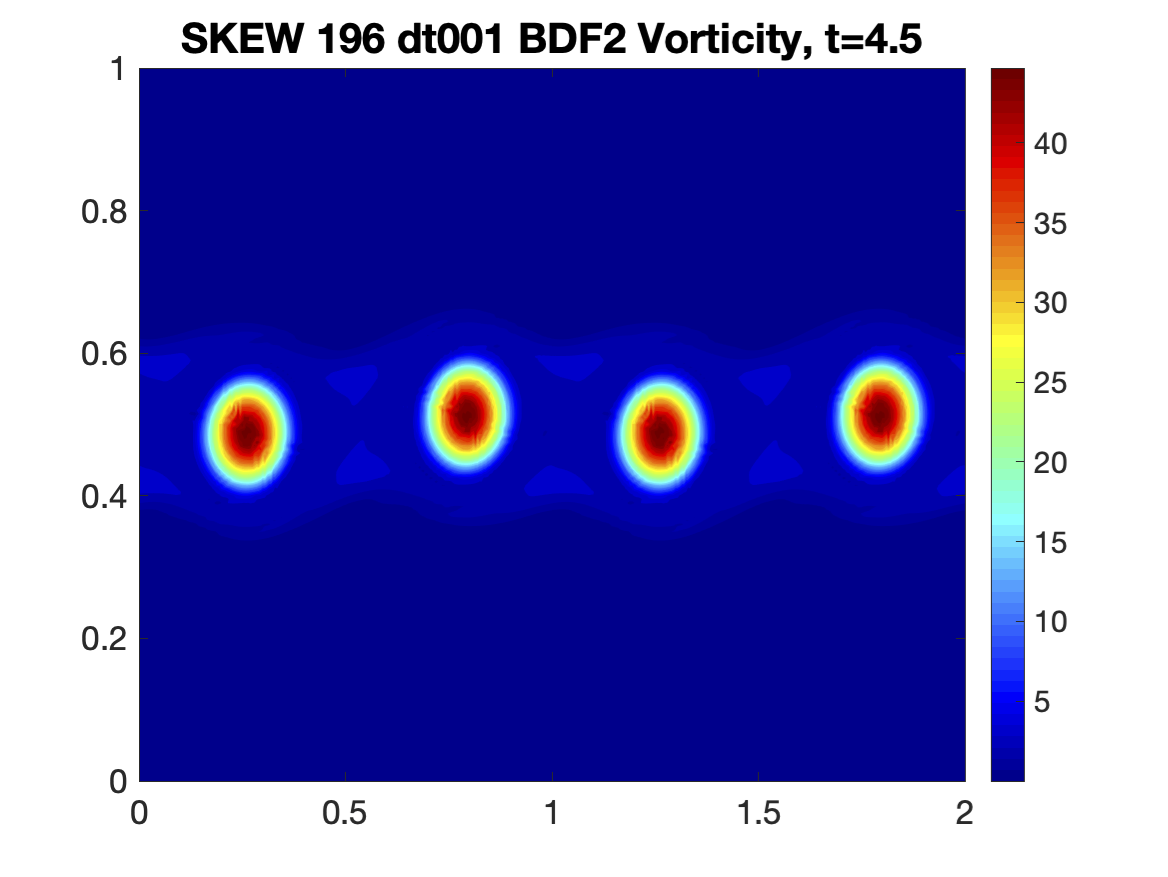} \\
\includegraphics[width=.4\textwidth, height=.18\textwidth,viewport=0 0 520 395, clip]{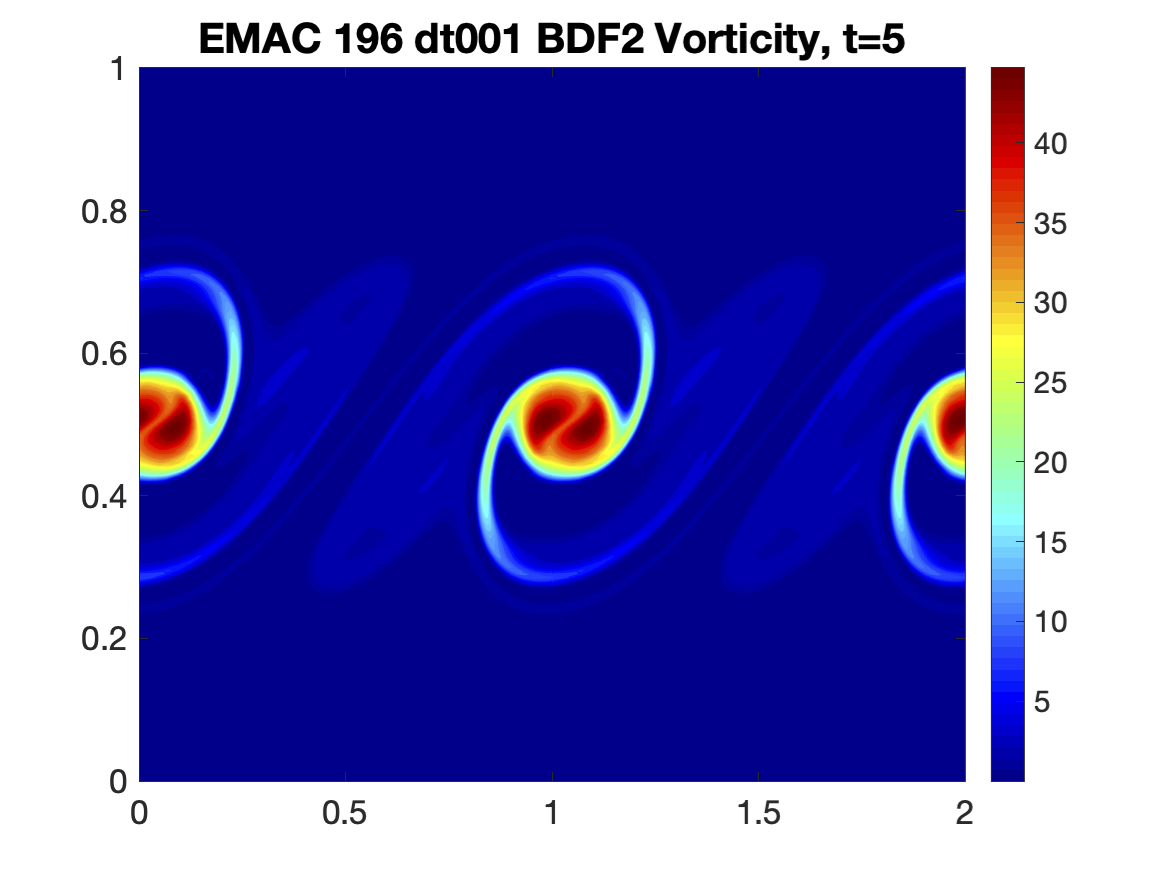}
\includegraphics[width=.4\textwidth, height=.18\textwidth,viewport=0 0 520 395, clip]{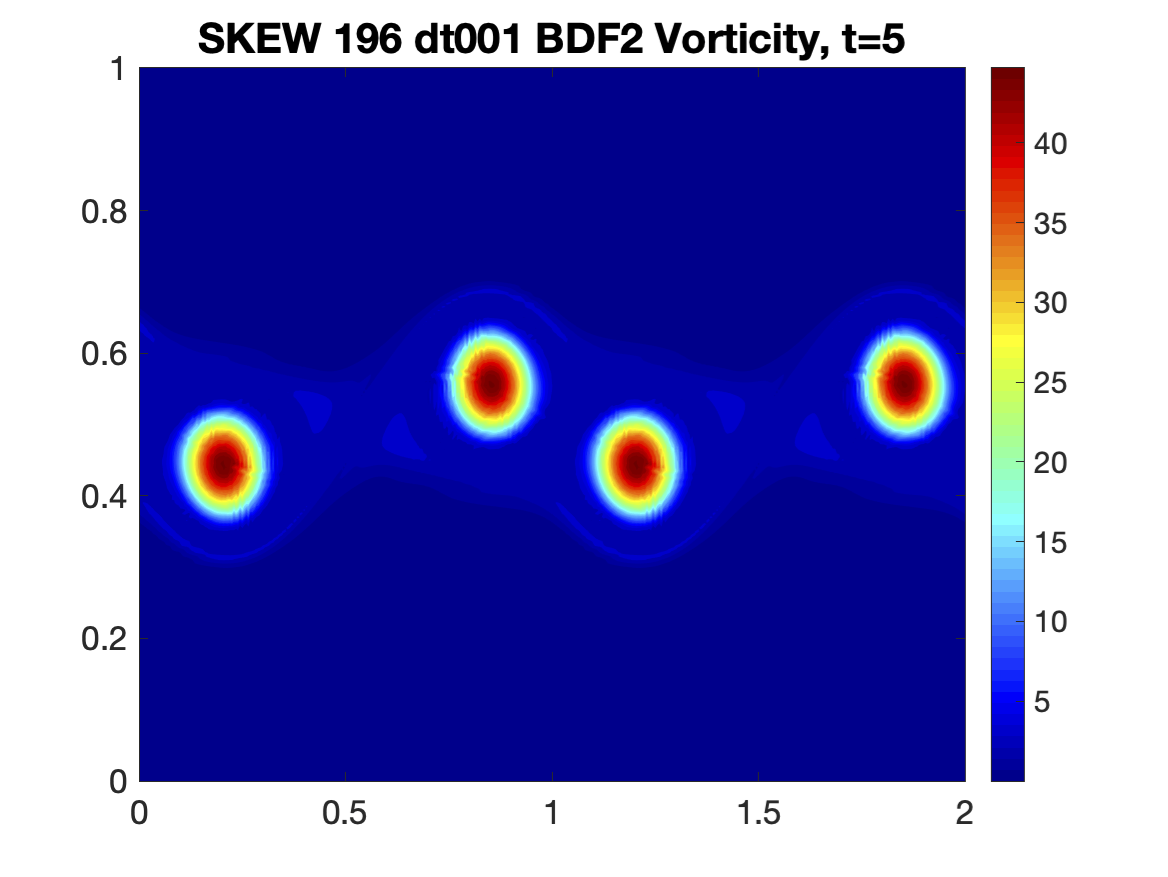} \\
\includegraphics[width=.4\textwidth, height=.18\textwidth,viewport=0 0 520 395, clip]{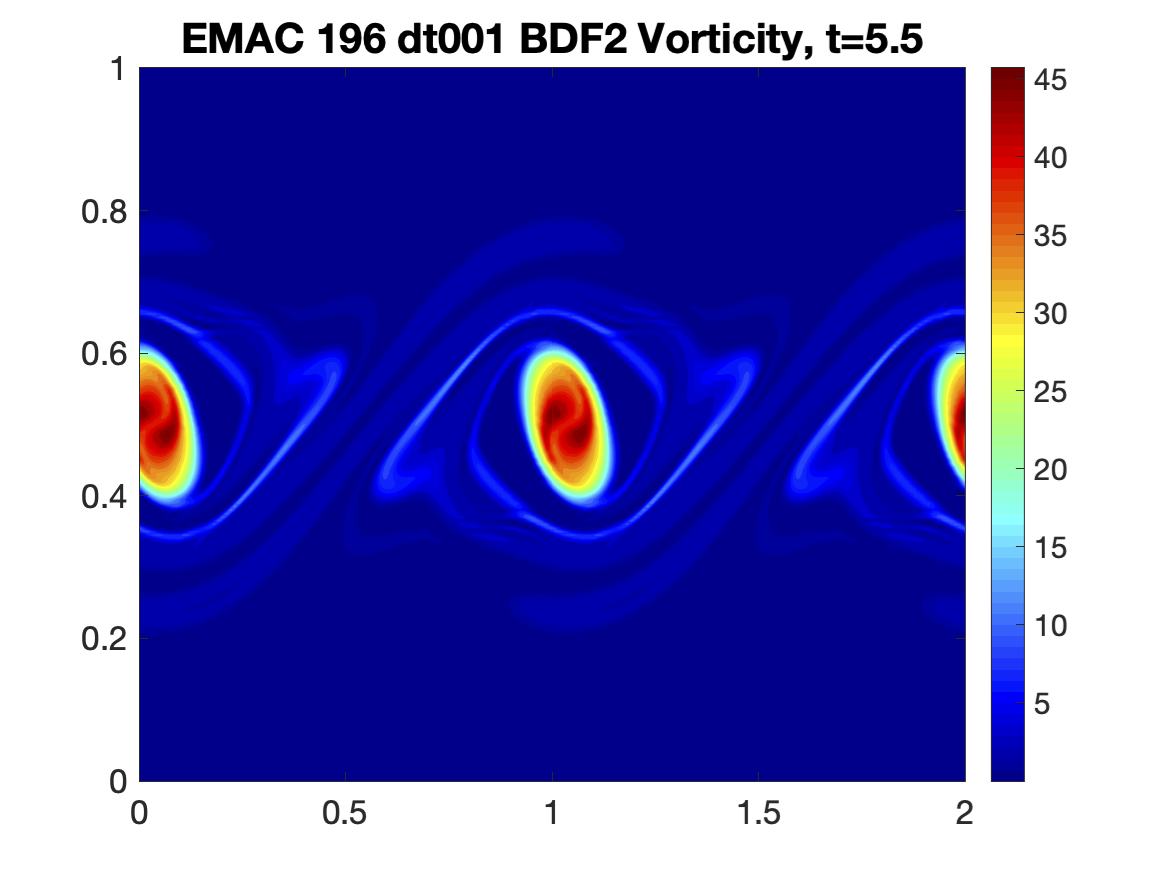}
\includegraphics[width=.4\textwidth, height=.18\textwidth,viewport=0 0 520 395, clip]{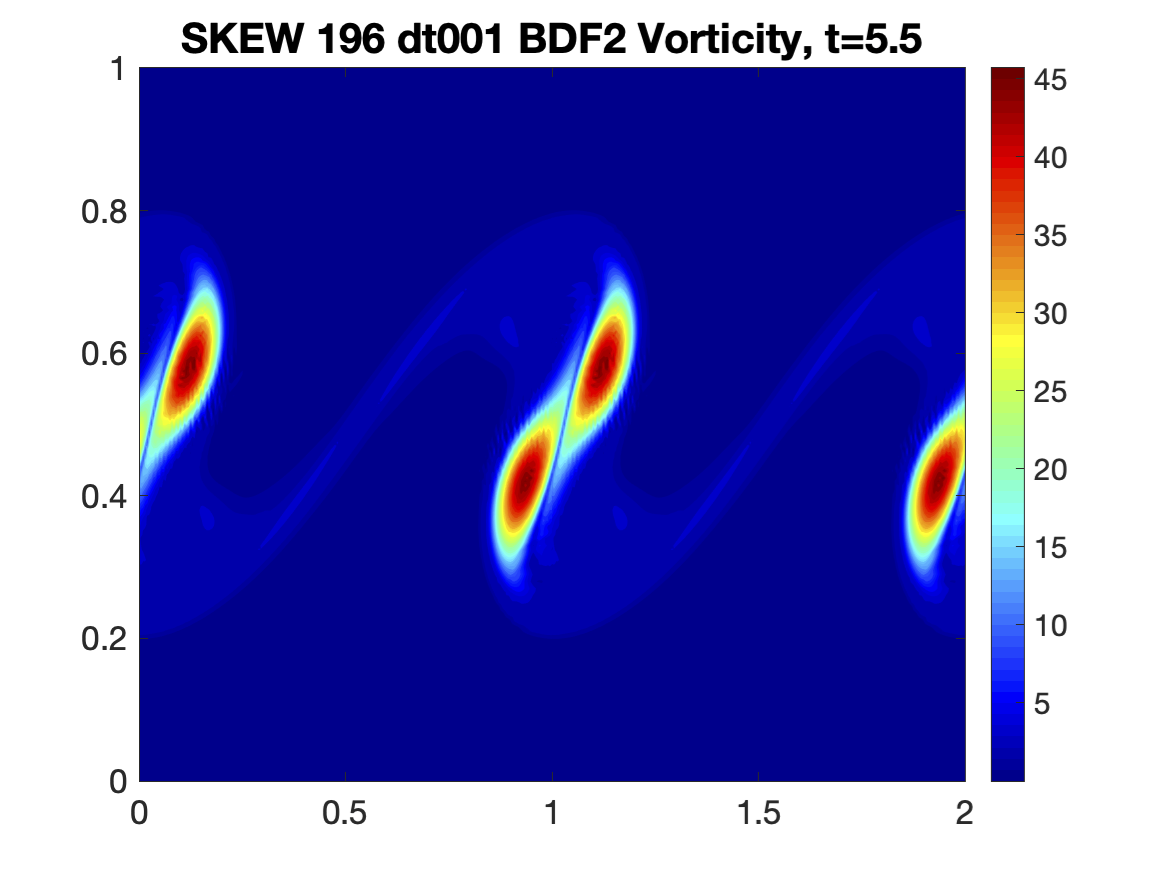} \\
\includegraphics[width=.4\textwidth, height=.18\textwidth,viewport=0 0 520 395, clip]{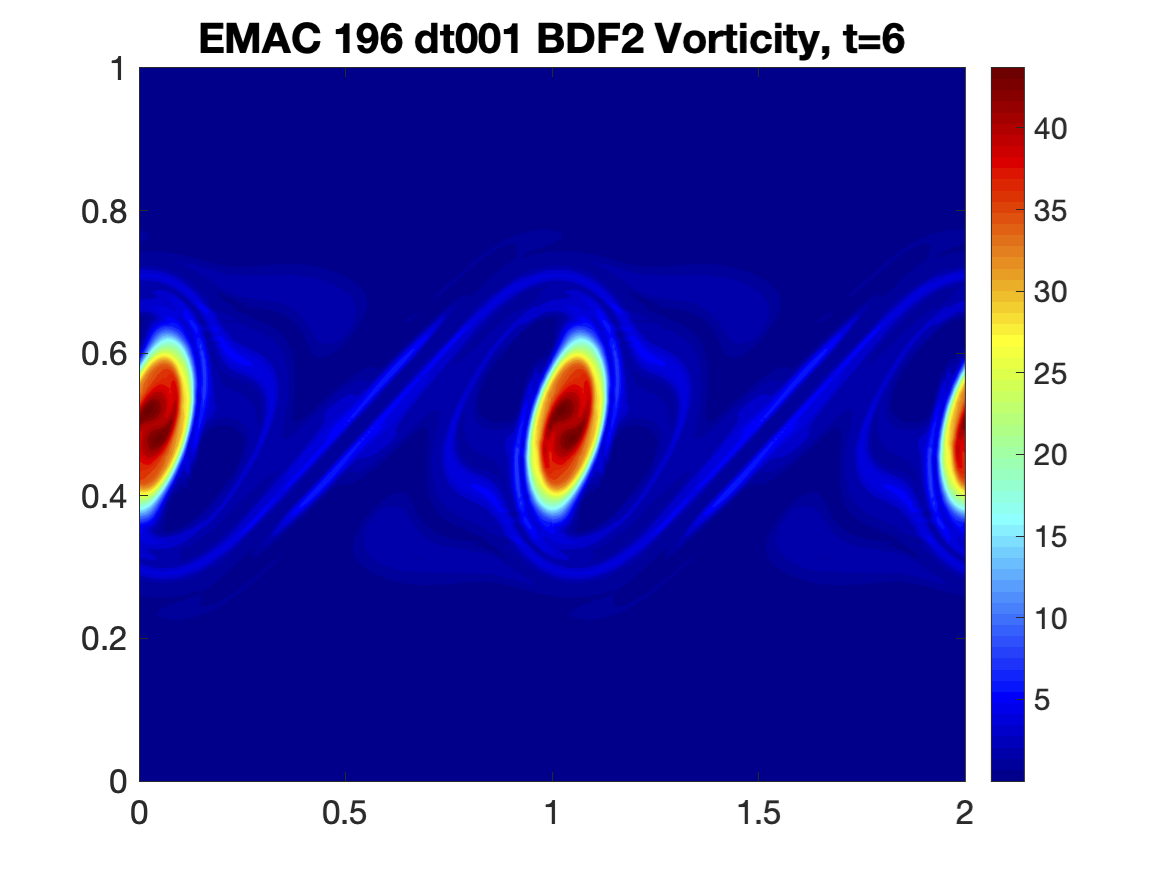}
\includegraphics[width=.4\textwidth, height=.18\textwidth,viewport=0 0 520 395, clip]{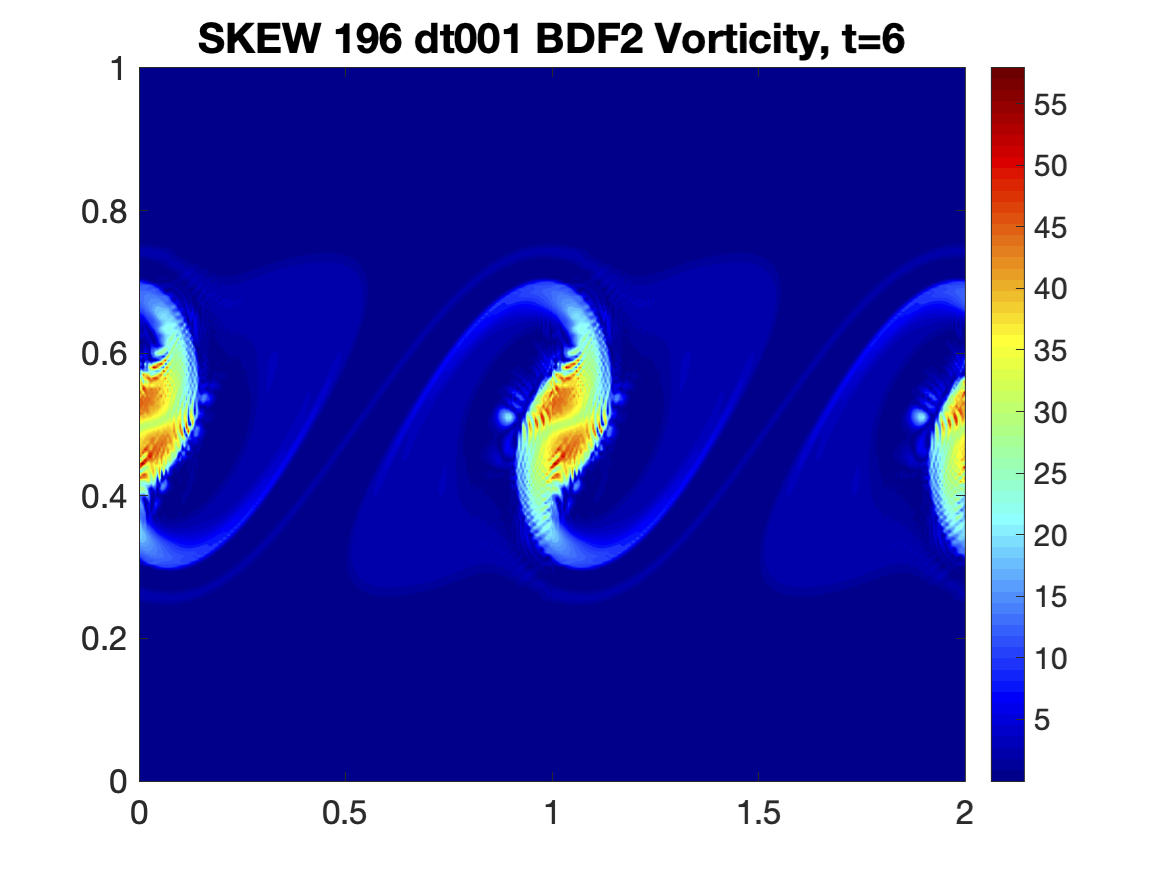} \\
\includegraphics[width=.4\textwidth, height=.18\textwidth,viewport=0 0 520 395, clip]{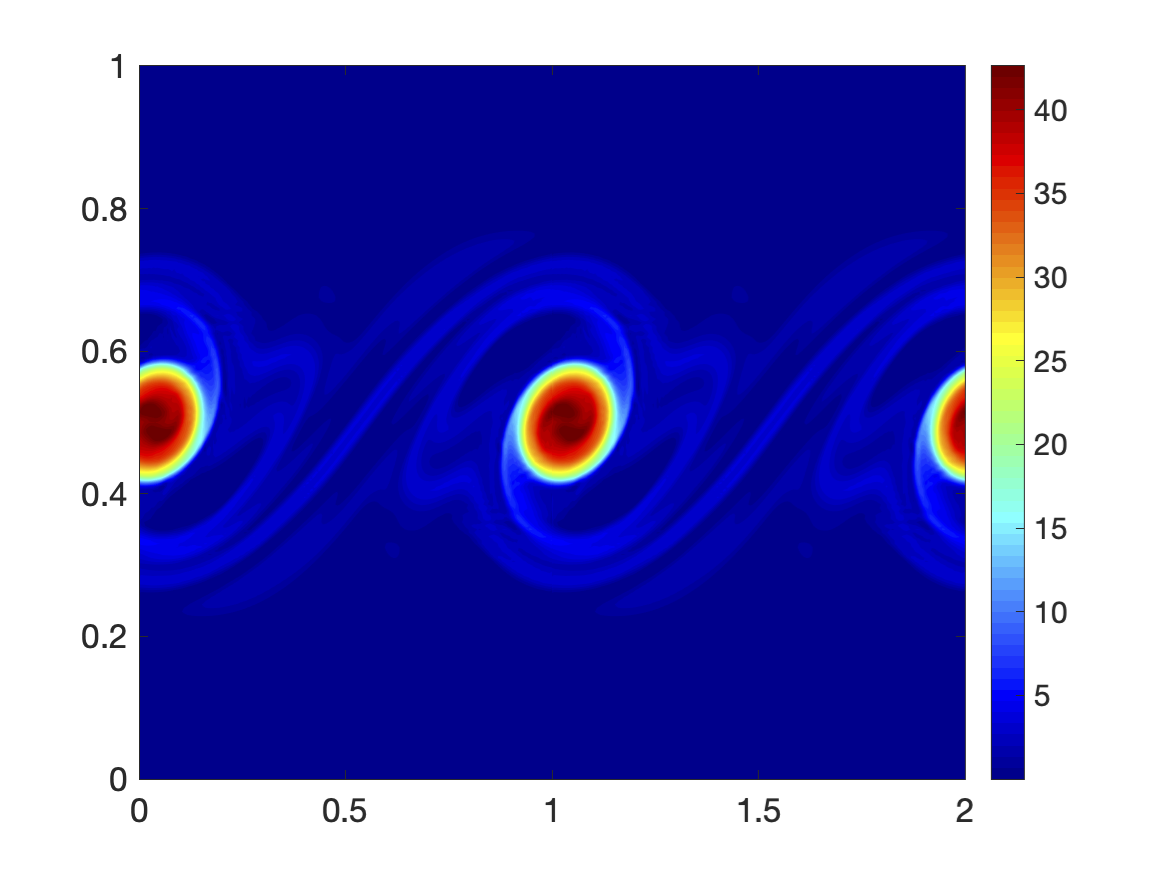}
\includegraphics[width=.4\textwidth, height=.18\textwidth,viewport=0 0 520 395, clip]{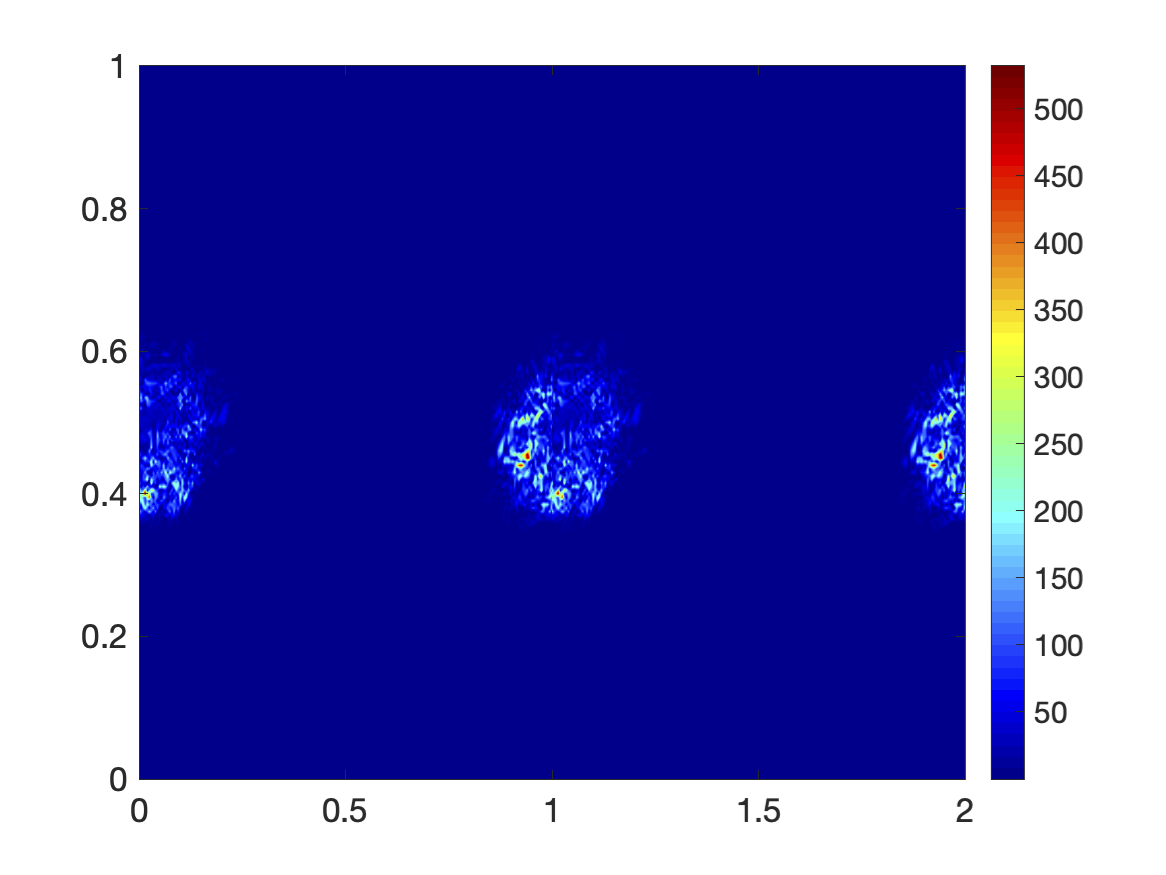} \\
\includegraphics[width=.4\textwidth, height=.18\textwidth,viewport=0 0 520 395, clip]{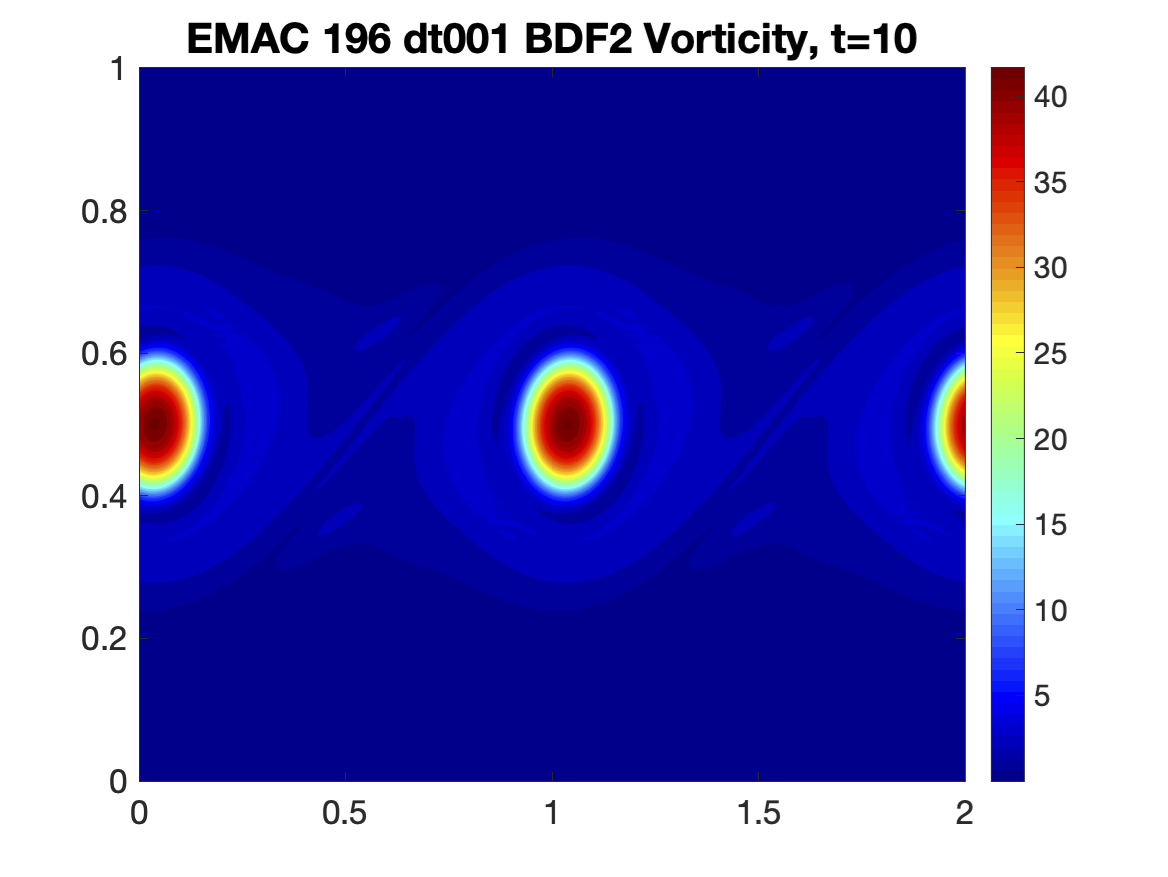}
\includegraphics[width=.4\textwidth, height=.18\textwidth,viewport=0 0 520 395, clip]{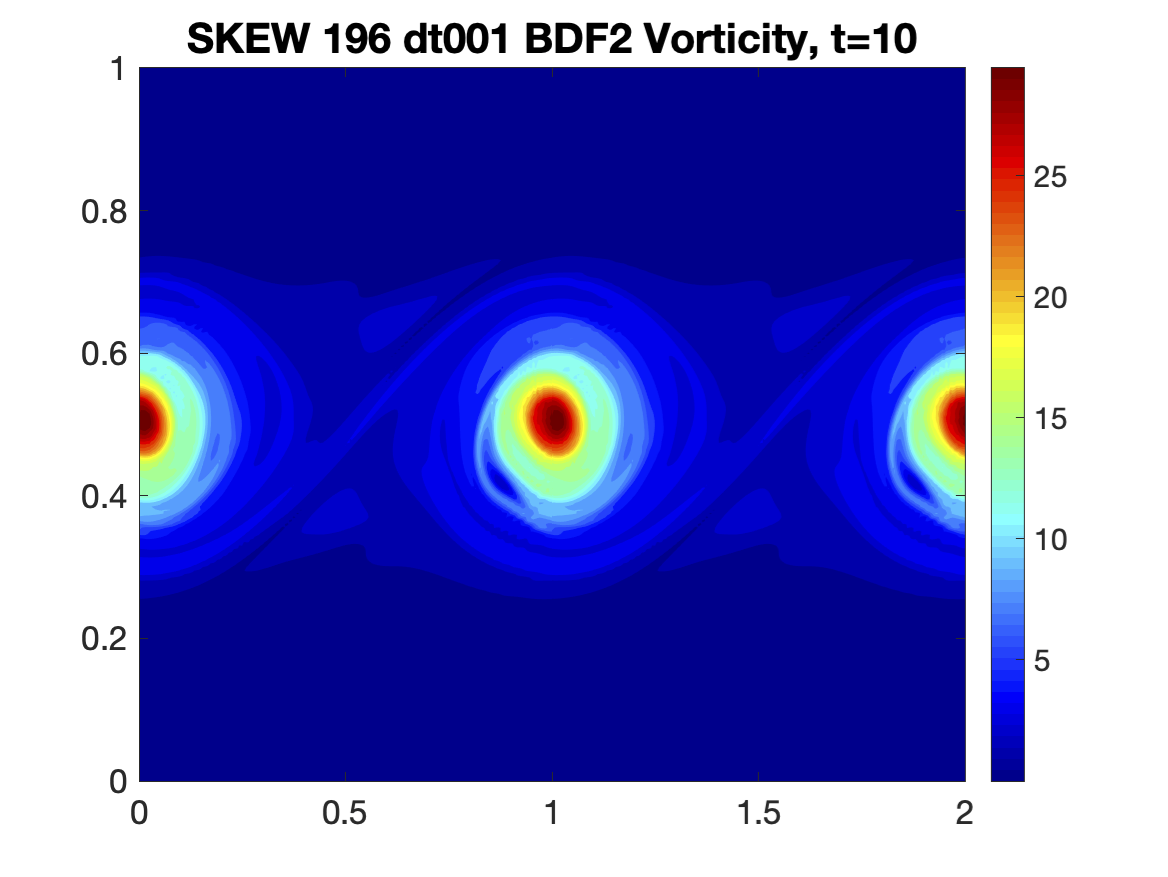}
\end{center}
\caption{\label{kh1000}
Shown above are absolute vorticity contours at for EMAC (left) and SKEW (right), at times t=1.0, 4.5, 5.0, 5.5, 6.0, 6.5 and 10.0 (from top to bottom), for $Re=1000$.}
\end{figure}

For our final test we consider a test problem from \cite{SJLLLS18} for 2D Kelvin-Helmholtz instability.  The domain is the unit square, with periodic
boundary conditions at $x=0,1$, representing an infinite extension in the horizontal direction.  At $y=0,1$, the no penetration boundary condition $\bu \cdot \bn=0$ is strongly enforced, along with a (natural) weak enforcement of the free-slip condition $(-\nu\nabla \bu \cdot\bn)\times \bn={\bf 0}$.  The initial condition is defined by
\[
\bu_0(x,y) = \left( \begin{array}{c} u_{\infty} \tanh\left( \frac{2y-1}{\delta_0} \right) \\ 0 \end{array} \right) + c_n \left( \begin{array}{c} \partial_y \psi(x,y) \\ -\partial_x \psi(x,y) \end{array} \right),
\]
where $\delta_0=\frac{1}{28}$ denotes the initial vorticity thickness, $u_{\infty}=1$ is a reference velocity, $c_n$ is a noise/scaling factor which is taken to be $10^{-3}$, and
\[
\psi(x,y) = u_{\infty} \exp \left( -\frac{(y-0.5)^2 }{\delta_0^2} \right) \left( \cos(8\pi x) + \cos(20\pi x) \right).
\]
The Reynolds number is defined by $Re=\frac{\delta_0 u_{\infty}}{\nu} = \frac{1}{28 \nu}$, and $\nu$ is defined by selecting $Re$.

We compute solutions for  $Re=1000$ using EMAC and SKEW formulations (we also tried CONV, but solutions became unstable and failed very early in the simulations).  Taylor-Hood $(P_2,P_1)$ elements are used without any stabilization for the spatial discretization, and BDF2 time stepping for the temporal discretization.  Solutions are computed  on a uniform triangulation with $h=\frac{1}{196}$ which provided 309K velocity degrees of freedom, up to end time $T=14$ using a time step size of $\Delta t=0.001$.   The nonlinear problems were resolved with Newton's method, and in most cases converged in 2 to 3 iterations.  For energy and enstrophy statistics, we compare results of SKEW and EMAC to results from Schroeder et al \cite{SJLLLS18} that used 5.5 million velocity degrees of freedom and a time step size of 3.6e-5.  We refer to results of \cite{SJLLLS18} as reference results.
Momentum and angular momentum reference statistics were not available in the literature for this test problem.

Plots of solutions are shown in figure \ref{kh1000} as enstrophy contours at different times.  It is known that the time evolution of Kelvin-Helmholtz flow is very sensitive \cite{SJLLLS18}, and so one should not read much into differences between the plots at particular times.  However, it is critical to notice that while EMAC appears to give resolved plots, the SKEW plots are not resolved and show oscillations.  This is most evident at $t=6.5$, but also is evident \MO{at $t=6$  and} also to a lesser extent at other times.  The SKEW solution plot at $t=6.5$ shows a clear problem in that a spike in enstrophy has (incorrectly) occurred.

Plots of energy, enstrophy, absolute total momentum, and angular momentum versus time are shown in figure \ref{conservation1000}.  Note that for enstrophy, figure \ref{conservation1000} shows two plots of the same data, with one zoomed in on the y axis. The EMAC energy matches the reference energy well for the entire simulation, and the EMAC enstrophy is accurate for $t<5$ and again for $t>10$.  EMAC enstrophy stays above the reference energy in $5<t<10$, seemingly because an inflection point for the reference solution curve occurs earlier than for EMAC; we note that a similar difference is seen in \cite{SJLLLS18} between this reference solution and solutions found on coarser meshes.  SKEW, on the other hand, is far from accurate, with the main symptom of the inaccuracy being a massive spike in the enstrophy between $t=6$ and $t=7$.  Not surprisingly, the energy prediction of SKEW also becomes inaccurate at $t=6$.  {\color{black} For momentum, EMAC preserves the zero momentum of the initial condition, while SKEW's momentum becomes more inaccurate with time and by $t=7$ reaches $10^{-3}$.  Angular momentum behavior for SKEW and EMAC are similar, with SKEW's appearing slightly delayed development between $t=4$ and $t=5$, which can be seen in  figure \ref{kh1000} at $t=4.5$, where EMAC has begun combining the eddies but SKEW has not.}

\begin{figure}[!h]
\begin{center}
\includegraphics[width=.32\textwidth, height=.25\textwidth,clip]{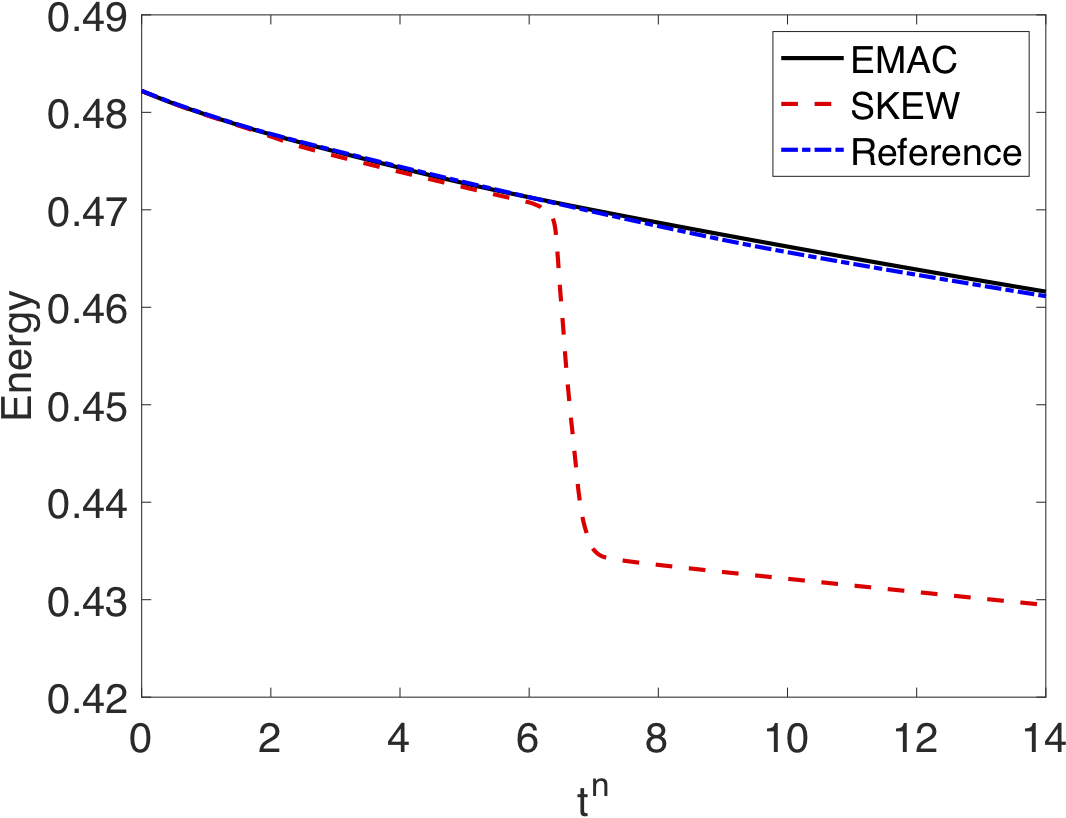} \
\includegraphics[width=.32\textwidth, height=.25\textwidth,clip]{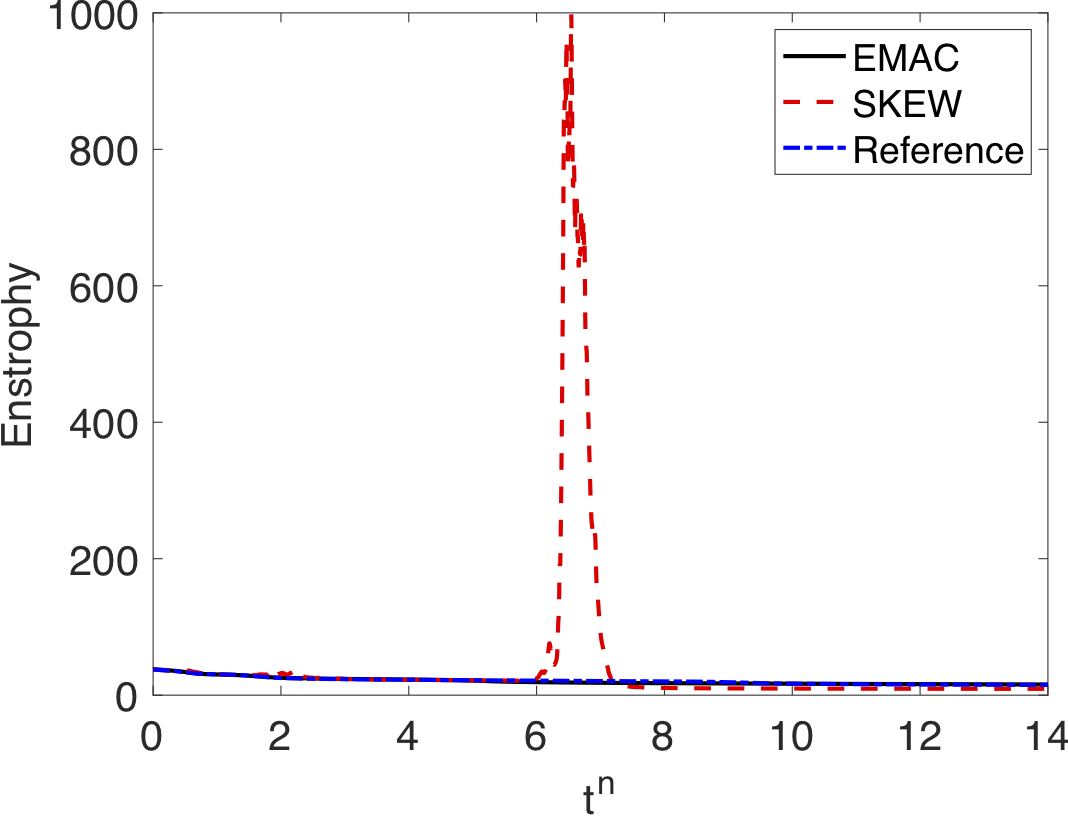} \
\includegraphics[width=.32\textwidth, height=.25\textwidth,clip]{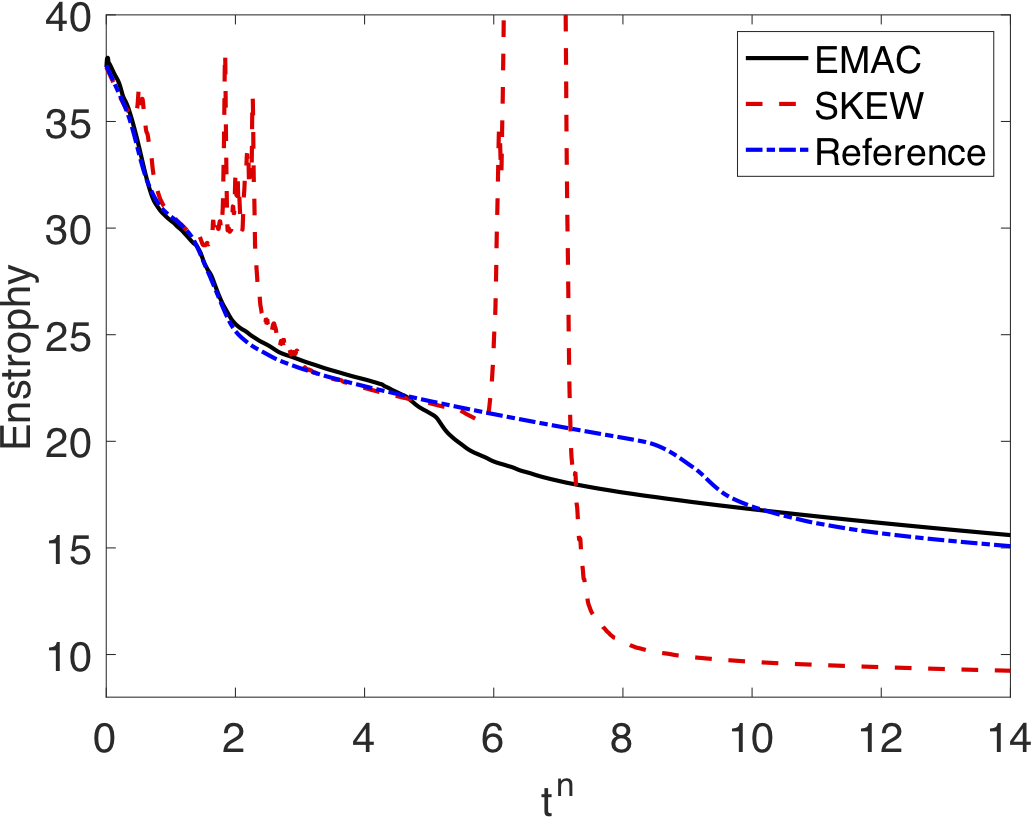} \\
\includegraphics[width=.32\textwidth, height=.25\textwidth,clip]{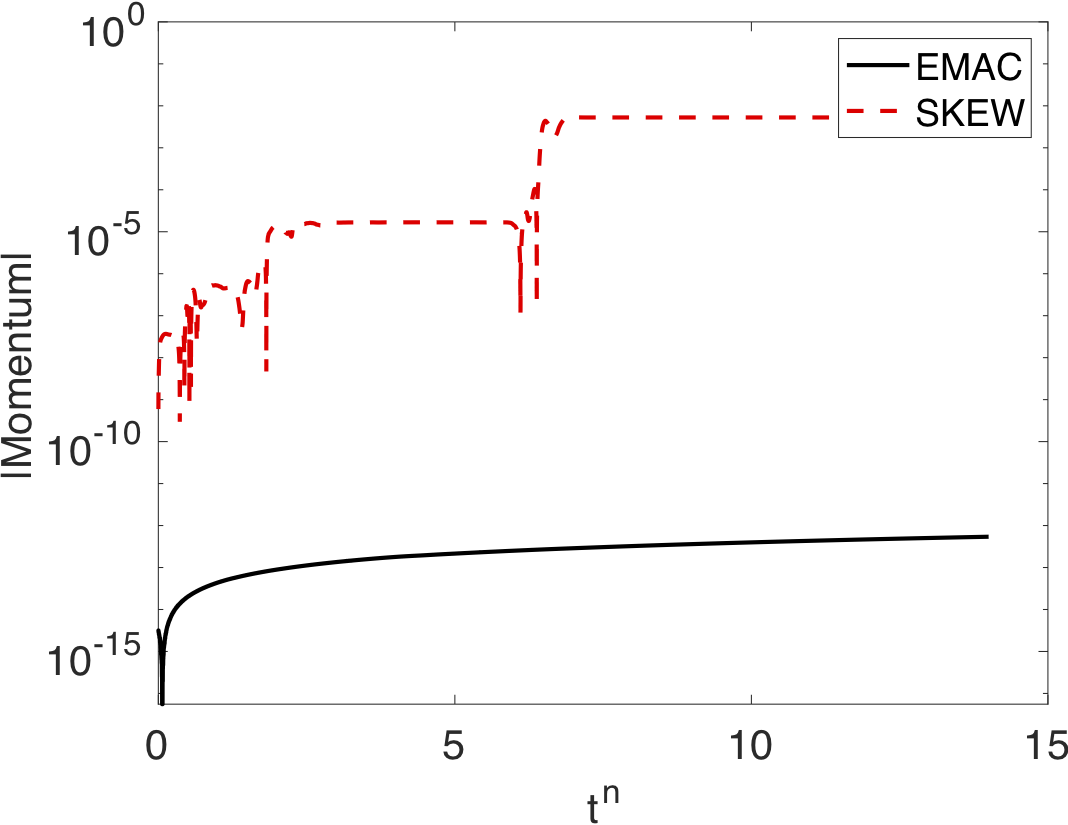} \
\includegraphics[width=.32\textwidth, height=.25\textwidth, clip]{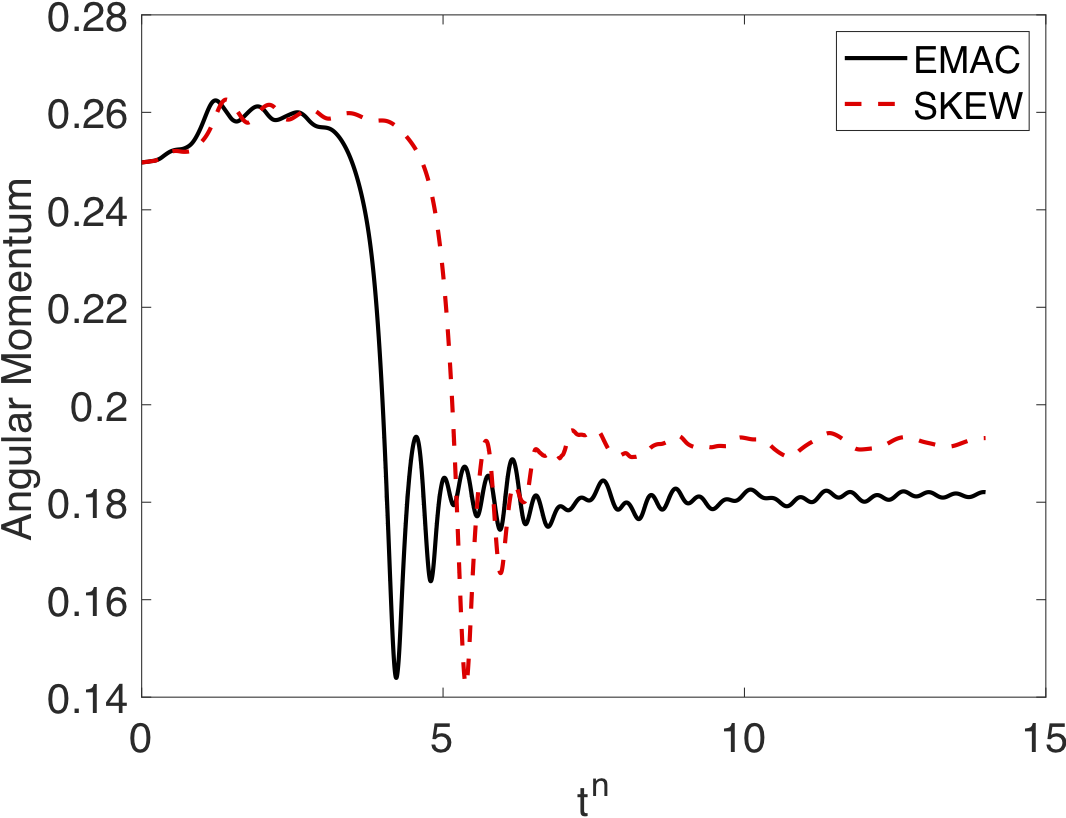}
\end{center}
\caption{\label{conservation1000}
Shown above are energy, enstrophy, and enstrophy zoomed in on y axis versus time, for $Re=1000$.}
\end{figure}

}

\section{Conclusions}

 We have given new analytical and numerical results that reveal more advantages of the EMAC
formulation in schemes for the incompressible Navier-Stokes equations. 
In particular, we have proven a better velocity error bound for EMAC which reduces (or removes) the dependence of the Gronwall constant on the Reynolds number compared to the classical result for SKEW, and we have also shown that an inaccurate momentum or angular momentum prediction (such as that of SKEW) creates a lower bound on the $L^2$ velocity error.  Both of these results suggest a better longer time accuracy of EMAC compared to SKEW, and the numerical results herein are in agreement.
 From a higher level, our results provide mathematical backing in this instance to the widely believed theory that more physically consistent schemes are in general more accurate over longer time intervals.  Future directions of this work may include extensions to incompressible multiphysics problems such as MHD and Boussinesq systems. More robust error bounds for other enhanced-physics methods through improved stability constants are also of interest.

\bibliographystyle{plain}

\end{document}